\documentclass[a4paper,12pt]{article}
\usepackage[cp1251]{inputenc}

\usepackage{epsfig}
\usepackage{amssymb}
\usepackage{amsthm}
\usepackage{graphicx}
\usepackage{amsfonts}
\usepackage{amsmath}
\usepackage{array}

\theoremstyle{plain}
\theoremstyle{plain}
\newtheorem{theorem}{Theorem}
\newtheorem*{theoremA}{Theorem A}
\newtheorem*{theoremB}{Theorem B}
\newtheorem*{theoremC}{Theorem C}

\newtheorem{pro}{Proposition}
\newtheorem{lemma}{Lemma}

\theoremstyle{definition}
\newtheorem{definition}{Definition}
\newtheorem*{notation}{Notation}
\newtheorem*{rem}{Remark}
\newtheorem*{ack}{Acknowledgements}

\newcommand{\bell}{\ensuremath{\boldsymbol{\ell}}}

\begin{document}

\author
{
R.\,R. Isangulov
\thanks {The research was partially supported by the~Russian Foundation for
Basic Research (grant 03--01--00104),
the~State Program of Support for Leading Scientific
Schools (grant 300.2003.1),
and INTAS (grant 03--51--3663).}
}

\title{Isospectral Flat $3$-Manifolds}

\date{}
\maketitle

\begin{abstract}
There is a well-known problem about isospectrality of Riemannian
manifolds: whether isospectral manifolds are isometric. In this work
we give an answer to this problem for compact flat 3-manifolds.
\end{abstract}


\section{Introduction}

\hspace*{\parindent}This paper deals with spectra of the
Laplace--Beltrami operator $\Delta=-\text{div grad}$ on compact
Riemannian manifold without boundary. Let us say that two
manifolds $M$ and $M'$ are {\it isospectral} if their spectra of
the Laplace--Beltrami operator on manifolds $M$ and $M'$ coincide.

There is a well-known problem about isospectrality of Riemannian
manifolds: whether isospectral manifolds are isometric.
In article \cite{Kac}, Kac framed this problem in the
question, "Can one hear the shape of a drum?" The answer is
negative in general. (See $\cite{Conway}$ for more information).

The first example of such manifolds was given by Milnor
\cite{Milnor} in 1964. He constructed a pair of 16-dimensional
isospectral non-isometric flat tori. In 1972
\mbox{McKean \cite{McKean}} showed that the cardinality of a set
of isospectral non-isometric compact Riemann surfaces is always
finite. In 1985 Sunada \cite{Sunada} discovered a general way of
constructing pairs of isospectral non-isometric Riemannian
manifolds. Using a new method in 1988 Brooks \cite{Brooks}
independently of Berger, Gauduchon and Mazet \cite{Berger} proved
that any two 2-dimensional flat tori are isospectral if and only
if they are isometric. In 1990 Schiemann \cite{Sch1} constructed a
pair of \mbox{4-dimensional} isospectral non-isometric flat tori.
In 1992 Conway and Sloane \cite{Conway&Sloane} constructed a simple
\mbox{4-parameter} family of pairs of isospectral 4-dimensional
lattices. Then in 1997 Schiemann \cite{Sch2} proved that any two
3-dimensional flat tori are isospectral if and only if they are
isometric. Recently, author \cite{Isan} independently of
Berger, Gauduchon and Mazet \cite{Berger} proved that any two flat
Klein bottles are isospectral if and only if they are isometric.
The main results and definitions concerning spectral theory of the
Laplace--Beltrami operator defined on Riemannian manifold can be
found, for example, in Buser \cite{Buser}.

In this paper we give an answer to this problem for
compact flat 3-manifolds. It is known \cite{Wolf} that there are
10 classes of pair-wise non-homeomorphic
compact flat 3-manifolds (6 classes are orientable and 4 classes
are non-orientable). Let $M_1,\ldots,M_6$ denote the 6 compact
orientable flat 3-manifolds and $N_1,\ldots,N_4$ denote the
non-orientable ones in the order given in \cite{Wolf}. To describe
the spectrum of manifold $M$ we use the trace function of the
manifold $M$ $\text{tr}(H_{M})=\int_{M}H_{M}(x,x,t)\,dM$, where
$H_{M}(x,y,t)$ is a fundamental solution of the heat equation on
$M$, $dM$ is the volume element. Using a method developed
in \cite{Isan} we have the following theorems.

\begin{theoremA}
The trace function of each compact flat $3$-manifolds $M_1,\ldots,
M_6,N_1,\ldots,N_4$ can be computed explicitly:
$\text{tr}(H_{M_i})=F_i$ for each manifold $M_i$, $i=1,\ldots,6$
and $\text{tr}(H_{N_j})=P_j$ for each manifold $N_j$,
$j=1,\ldots,4$, where the functions $F_i$ and $P_j$ are given in
section $6$, Table $1$.
\end{theoremA} 

\begin{theoremB}
Any two homeomorphic compact flat $3$-manifolds are isospectral if
and only if they are isometric.
\end{theoremB}\vspace{0.2cm}

\begin{theoremC}
There is a unique family of pairs of isospectral non-homeomorphic flat
$3$-ma\-ni\-folds
which consists of manifolds $M_4$ and $M_6$.
\end{theoremC}\vspace{0.2cm}

\begin{rem}
In order to prove Theorem A in case of $M_1$ (3-dimensional flat
torus) we use an explicit formula for the trace function of $M_1$
obtained in \cite{Brooks}.
To prove Theorem B in case of $M_1$ we
use a result obtained in \cite{Sch2}. All other cases are proved
using completely different methods.
\end{rem}\vspace{0.2cm}

In order to prove Theorem A we use a relationship between the
fundamental solutions of the heat equation on manifold $M$ and
on regular covering of this manifold $M$ (Lemma \ref{Brooks}).

To prove Theorem B we use Theorem A and the following statements:

(i) Spectrum of the manifold $M$ determines exactly the trace
function $\text{tr}(H_{M})$ and, conversely, spectrum of the
manifold $M$ can be determined by the trace function
$\text{tr}(H_{M})$ (Proposition \ref{pro1} and Proposition \ref{pro2}).

(ii) The compact flat 3-manifold $M$ can be determined up to
isometry by the trace function
$\text{tr}(H_{M})$ (Lemma \ref{mainlemma}).

Theorem C follows immediately from Theorem A.

Throughout this paper we use \cite{Buser} and \cite{Wolf} for
standard references. Section 2 contains preliminaries. Section 3,
Section 4 and Section 5 contain the proof of Theorem~A, Theorem~B and Theorem~C
correspondingly. Section 6 contains the table of trace functions $\text{tr}(H_M)$
and the fundamental sets of compact flat $3$-manifolds.
Results of this work are the subject of master thesis
"Isospectral flat $3$-manifolds" defended by the author at Novosibirsk
state university in 2002. However, as this work was being written, we found
that similar results were obtained independently in recent preprints
\cite{Conway&Rossetti} and \cite{Rossetti&Conway}.
\vspace{0.1cm}

\begin{ack}
I wish to thank my scientific supervisor Prof. Alexander D.
Mednykh for introducing this problem to me and for continuous
encouragement.
\end{ack}


\section{Preliminaries}\label{prelim}

In this section we will briefly summarize some of the basic facts
about compact flat \mbox{3-manifolds} (see \cite{Wolf} for details) and
about spectral theory of the Laplace--Beltrami operator
\mbox{$\Delta=-\text{div grad}$}
defined on compact Riemannian manifold without boundary (see \cite{Buser} for
details).

\subsection{Compact flat $3$-manifolds}

Let $E(n)$ denote the group of rigid motions of $\mathbb{R}^{n}$.
Every rigid motion consists of a translation $t_{a}$ by a vector
$a$ followed by a rotation $A$. Write the motion $(A,t_{a})$.
Clearly $A$ is an element of $O(n)$ and $a$ is an arbitrary vector
in $\mathbb{R}^{n}$. Thus the Euclidean group $E(n)$ is the
semi--direct product of $O(n)$ and $\mathbb{R}^{n}$ satisfying the
following product rule:
$$(A,t_a)(B,t_b)=(AB,t_{Ab+a}).$$

\begin{definition}
A {\it flat compact connected $n$-manifold} $M^{n}$ is the orbit space
of $\mathbb{R}^{n}$ by the fixed point free properly discontinuous action of a
discrete subgroup $\Gamma \in E(n)$, $M^n=\mathbb{R}^{n}/\Gamma$.
It admits a covering by the torus
$T^n=\mathbb{R}^{n}/\Gamma^{\ast}$, where $\Gamma^{\ast}$ is a normal subgroup
of rank $n$ of a finite index, $\Gamma^{\ast} < \Gamma$. Moreover,
one can choose $\Gamma^{\ast}=\Gamma\cap\mathbb{R}^{n}$. We note
also (see \cite[chapter~3]{Wolf}) that $\Gamma$ has no
non-trivial element of finite order.
\end{definition}

\begin{definition}
The group of deck
transformations $\Psi$ in the covering $T^n \rightarrow M^n$ is called
the {\it holonomy group} of $M^n$,
$\Psi=\Gamma/\Gamma^\ast$.
\end{definition}

The following results can be found in \cite{Wolf}.

\begin{theorem}\label{orientable}
There are just $6$ affine diffeomorphism classes of compact
connected orientable flat $3$-manifolds. They are represented by the
manifolds $\mathbb{R}^{3}/\Gamma$ where $\Gamma$ is one of the six
groups $M_i$ given below. Here $\Lambda$ is the translation
lattice, $\{ a_1,a_2,a_3 \}$ are its generators, $t_i=t_{a_i}$,
and $\Psi=\Gamma/\Gamma^\ast$ is the holonomy. \vspace{0.1cm}

\noindent $M_1$. \hspace{3pt}$\Psi=\{1\}$ and $\Gamma$ is
generated by the translations $\{ t_1,t_2,t_3\}$ with $\{ a_i\}$
linearly independent.\vspace{0.1cm}

\noindent $M_2$. \hspace{3pt}$\Psi=\mathbb{Z}_2$ and $\Gamma$ is
generated by $\{\alpha,t_1,t_2,t_3\}$ where $\alpha^2=t_1$,
$\alpha t_2 \alpha^{-1}=t_2^{-1}$ and

$\alpha t_3 \alpha^{-1}=t_3^{-1}$; $a_1$ is orthogonal to $a_2$
and $a_3$ while $\alpha=(A,t_{a_1/2})$ with $A(a_1)=a_1$,

$A(a_2)=-a_2$, $A(a_3)=-a_3$.\vspace{0.1cm}

\noindent $M_3$. \hspace{3pt}$\Psi=\mathbb{Z}_3$ and $\Gamma$ is
generated by $\{\alpha,t_1,t_2,t_3\}$ where $\alpha^3=t_1$,
$\alpha t_2 \alpha^{-1}=t_3$ and

$\alpha t_3 \alpha^{-1}=t_2^{-1}t_3^{-1}$; $a_1$ is orthogonal to
$a_2$ and $a_3$, $||a_2||=||a_3||$ and $\{a_2, a_3 \}$ is a
hexagonal

plane lattice, and $\alpha=(A,t_{a_1/3})$ with $A(a_1)=a_1$,
$A(a_2)=a_3$, $A(a_3)=-a_2-a_3$.\vspace{0.1cm}

\noindent $M_4$. \hspace{3pt}$\Psi=\mathbb{Z}_4$ and $\Gamma$ is
generated by $\{\alpha,t_1,t_2,t_3\}$ where $\alpha^4=t_1$,
$\alpha t_2 \alpha^{-1}=t_3$ and

$\alpha t_3 \alpha^{-1}=t_2^{-1}$;
$a_i$ are mutually orthogonal with $||a_2||=||a_3||$ while
$\alpha=(A,t_{a_1/4})$ with

$A(a_1)=a_1$,
$A(a_2)=a_3$, $A(a_3)=-a_2$.\vspace{0.1cm}

\noindent $M_5$. \hspace{3pt}$\Psi=\mathbb{Z}_6$ and $\Gamma$ is
generated by $\{\alpha,t_1,t_2,t_3\}$ where $\alpha^6=t_1$,
$\alpha t_2 \alpha^{-1}=t_3$ and

$\alpha t_3
\alpha^{-1}=t_2^{-1}t_3$;
$a_1$ is orthogonal to $a_2$ and $a_3$, $||a_2||=||a_3||$ and
$\{a_2, a_3 \}$ is a hexagonal

plane lattice, and
$\alpha=(A,t_{a_1/6})$ with $A(a_1)=a_1$, $A(a_2)=a_3$,
$A(a_3)=-a_3-a_2$.\vspace{0.1cm}

\noindent $M_6$. \hspace{3pt}$\Psi=\mathbb{Z}_2\times\mathbb{Z}_2$
and $\Gamma$ is generated by $\{\alpha,\beta,\gamma,t_1,t_2,t_3\}$
where $\gamma\beta\alpha=t_1 t_3$ and $\alpha^2=t_1$,

$\alpha t_2
\alpha^{-1}=t_2^{-1}$, $\alpha t_3 \alpha^{-1}=t_3^{-1}$; $\beta
t_1 \beta^{-1}=t_1^{-1}$, $\beta^2=t_2$, $\beta t_3
\beta^{-1}=t_3^{-1}$; $\gamma t_1 \gamma^{-1}=t_1^{-1}$,

$\gamma
t_2 \gamma^{-1}=t_2^{-1}$, $\gamma^2=t_3$; the $a_i$ are mutually
orthogonal and

$\alpha=(A,t_{a_1/2})$ with $A(a_1)=a_1$, $A(a_2)=-a_2$,
$A(a_3)=-a_3$;

$\beta=(B,t_{(a_2+a_3)/2})$ with $B(a_1)=-a_1$, $B(a_2)=a_2$,
$B(a_3)=-a_3$;

$\gamma=(C,t_{(a_1+a_2+a_3)/2})$ with $C(a_1)=-a_1$,
$C(a_2)=-a_2$, $C(a_3)=a_3$.

\end{theorem}
\vspace{0.2cm}

\begin{theorem}\label{nonorientable}
There are just $4$ affine diffeomorphism classes of compact
connected non-orientable flat $3$-manifolds. They are represented by
the manifolds $\mathbb{R}^{3}/\Gamma$ where $\Gamma$ is one of the
$4$ groups $N_i$ given below. Here $\Lambda$ is the translation
lattice, $\{ a_1,a_2,a_3 \}$ are its generators, $t_i=t_{a_i}$,
$\Psi=\Gamma/\Gamma^\ast$ is the holonomy, and
$\Gamma_0=\Gamma\cap \text{SO}(3)\cdot\mathbb{R}^3$ so that
$\mathbb{R}^3/\Gamma_0 \rightarrow \mathbb{R}^3/\Gamma$ is the
$2$-fold orientable Riemannian covering.\vspace{0.1cm}

\noindent $N_1$. \hspace{3pt}$\Psi=\mathbb{Z}_2$ and $\Gamma$ is
generated by $\{\varepsilon,t_1,t_2,t_3\}$ where
$\varepsilon^2=t_1$, $\varepsilon t_2 \varepsilon^{-1}=t_2$ and
$\varepsilon t_3 \varepsilon^{-1}=t_3^{-1}$;

$a_1$ and $a_2$ are orthogonal to $a_3$ while
$\varepsilon=(E,t_{a_1/2})$ with $E(a_1)=a_1$, $E(a_2)=a_2$,

$E(a_3)=-a_3$. $\Gamma_0$ is generated by
$\{t_1,t_2,t_3\}$.\vspace{0.1cm}

\noindent $N_2$. \hspace{3pt}$\Psi=\mathbb{Z}_2$ and $\Gamma$ is
generated by $\{\varepsilon,t_1,t_2,t_3\}$ where
$\varepsilon^2=t_1$, $\varepsilon t_2 \varepsilon^{-1}=t_2$ and

$\varepsilon t_3 \varepsilon^{-1}=t_1 t_2 t_3^{-1}$;
the orthogonal projection of $a_3$ on the $(a_1,a_2)$--plane is
$(a_1+a_2)/2$,

while $\varepsilon=(E,t_{a_1/2})$
with $E(a_1)=a_1$, $E(a_2)=a_2$, $E(a_3)=a_1+a_2-a_3$.

$\Gamma_0$ is generated by $\{t_1,t_2,t_3\}$.\vspace{0.1cm}

\noindent $N_3$. \hspace{3pt}$\Psi=\mathbb{Z}_2\times\mathbb{Z}_2$
and $\Gamma$ is generated by $\{\varepsilon,\alpha,t_1,t_2,t_3\}$
where $\alpha^2=t_1$, $\varepsilon^2=t_2$,
$\varepsilon\alpha\varepsilon^{-1}=t_2\alpha$,

$\alpha t_2 \alpha^{-1}=t_2^{-1}$, $\alpha t_3
\alpha^{-1}=t_3^{-1}$, $\varepsilon t_1 \varepsilon^{-1}=t_1$ and
$\varepsilon t_3 \varepsilon^{-1}=t_3^{-1}$; the $a_i$ are
mutually

orthogonal and

$\alpha=(A,t_{a_1/2})$ with $A(a_1)=a_1$, $A(a_2)=-a_2$,
$A(a_3)=-a_3$;

$\varepsilon=(E,t_{a_2/2})$ with $E(a_1)=a_1$, $E(a_2)=a_2$,
$E(a_3)=-a_3$.

$\Gamma_0$ is generated by $\{\alpha,t_1,t_2,t_3\}$.
\vspace{0.1cm}

\noindent $N_4$. \hspace{3pt}$\Psi=\mathbb{Z}_2\times\mathbb{Z}_2$
and $\Gamma$ is generated by $\{\varepsilon,\alpha,t_1,t_2,t_3\}$
where $\alpha^2=t_1$, $\varepsilon^2=t_2$,
$\varepsilon\alpha\varepsilon^{-1}=t_2t_3\alpha$,

$\alpha t_2 \alpha^{-1}=t_2^{-1}$, $\alpha t_3
\alpha^{-1}=t_3^{-1}$, $\varepsilon t_1 \varepsilon^{-1}=t_1$ and
$\varepsilon t_3 \varepsilon^{-1}=t_3^{-1}$; the $a_i$ are
mutually

orthogonal and

$\alpha=(A,t_{a_1/2})$ with $A(a_1)=a_1$, $A(a_2)=-a_2$,
$A(a_3)=-a_3$;

$\varepsilon=(E,t_{(a_2+a_3)/2})$ with $E(a_1)=a_1$, $E(a_2)=a_2$,
$E(a_3)=-a_3$.

$\Gamma_0$ is generated by $\{\alpha,t_1,t_2,t_3\}$.

\end{theorem}
\vspace{0.2cm}

\begin{theorem}\label{isometry}
Let $M$ be a compact flat $3$-manifold, $M=\mathbb{R}^3/\Gamma$.
Then the following statements are hold.\vspace{0.1cm}

\rm (i) \it If $M$ is a torus $M_1$ then all the possible isometry classes
are parameterized by elements from
$\text{SL}(3,\mathbb{Z})\backslash \text{GL}^{+}(3,\mathbb{R})/\text{SO}(3).$
\vspace{0.1cm}

\rm (ii) \it
If $M$ is a manifold $M_2$ then
isometry classes can be parameterized by
$||a_1||$ and $2$-dimensional lattice $\{a_2,\,a_3\}$ generated by
vectors $a_2$ and $a_3$.
\vspace{0.1cm}

\rm (iii) \it If $M$ is a manifold $M_3$, $M_4$ and $M_5$ then isometry classes
can be parameterized by $||a_1||$ and $||a_3||$.
\vspace{0.1cm}

\rm (iv) \it If $M$ is a manifold $M_6$ then the lattice
$\Lambda=\{a_1,a_2,a_3\}$ is a right-angled one and all the
possible isometry classes are parameterized by unordered set of
three $(||a_1||,||a_2||,||a_3||)$.\vspace{0.1cm}

\rm (v) \it If $M$ is a manifold $N_1, N_2$ then isometry classes are
parameterized by pair $(p,||a_3||)$, where
$p$
is an element corresponding to $2$-dimensional lattice
$\{a_1,\,a_2\}$ generated by vectors $a_1$ and $a_2$.
\vspace{0.1cm}

\rm (vi) \it If $M$ is a manifold $N_3, N_4$ then the lattice
$\Lambda=\{a_1,a_2,a_3\}$ is a right-angled one and all the
possible cases are parameterized by ordered set of three
$(||a_1||,||a_2||,||a_3||)$.

\end{theorem}
\vspace{0.2cm}

It was shown by Luft and Sjerve in \cite{Luft} that compact flat
3-manifolds admit the following hierarchy of (topological)
coverings (see Fig. 1), where number next to arrow is the number
of folds in the corresponding regular covering. We need the following
well-known observation.

\begin{rem}
Let $\pi: M' \rightarrow M$ be a
topological covering and $M$ is endowed with flat metric then
this metric can be lifted to the flat metric on $M'$ in such a
way that $\pi$ becomes a local isometry. The inverse is not true:
flat metic on $M'$ is not necessary a lifting of flat metric on
$M$ via $\pi$.
\end{rem}

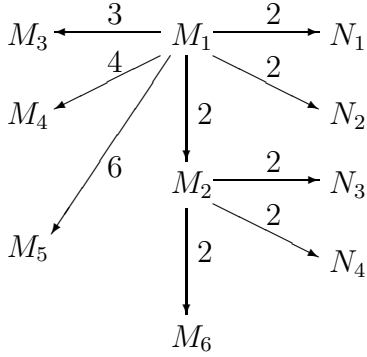
\begin{figure}
\centering
\begin{picture}(200,200)
\put(100,160){$M_1$}

\put(116,165){\vector(1,0){40}} \put(160,160){$N_1$}
\put(136,168){2}

\put(96,165){\vector(-1,0){40}} \put(38,160){$M_3$}
\put(76,168){3}

\put(116,156){\vector(2,-1){40}}\put(160,130){$N_2$}
\put(136,148){2}

\put(96,156){\vector(-2,-1){40}}\put(38,130){$M_4$}
\put(76,150){4}

\put(106,156){\vector(0,-1){40}} \put(100,104){$M_2$}
\put(110,130){2}

\put(116,109){\vector(1,0){40}} \put(160,104){$N_3$}
\put(136,112){2}

\put(116,100){\vector(2,-1){40}}\put(160,74){$N_4$}
\put(136,92){2}

\put(106,98){\vector(0,-1){40}} \put(100,46){$M_6$}
\put(110,78){2}

\put(100,156){\vector(-2,-3){45}}\put(38,80){$M_5$}
\put(76,110){6}

\end{picture}
\caption{Hierarchy of (topological) coverings}

\end{figure}
\vspace{0.2cm}

\subsection{Spectral theory of the Laplace--Beltrami operator}

In this section we will briefly give the main results concerning
spectral theory of the Laplace--Beltrami operator defined on
Riemannian manifold. It can be found in \cite[chapter 7]{Buser}.

\begin{definition}
{\it Laplace--Beltrami operator} (or {\it
laplacian} for short) is an operator $$\Delta(f)=-\text{div grad}f.$$
\end{definition}

Let $M$ be an arbitrary complete $m$-dimensional Riemannian
manifold. The Laplace--Beltrami operator has the coordinate expression

$$\Delta u=-\frac{1}{\sqrt{g}} \sum_{j,k=1}^{m}
\partial_{j}(g^{jk}\sqrt{g}\,\partial_{k} u),\hspace{10pt} u\in
C^{\infty}(M),$$

\noindent where $\partial_{j}$ denotes the partial differentiation
with respect to the $j$-th coordinate; $g=\text{det}(g_{ij})$
is the determinant of the metric tensor, where the $g_{ij}$
are the components of the metric tensor with respect to the local
coordinates and the $g^{ij}$ are the components of the inverse tensor.

The function space in question is the Hilbert space $L_{2}(M)$ of
square integrable functions
\mbox{$f:M\rightarrow \mathbb{R}$} on the
compact connected unbordered Riemannian manifold $M$ with inner
product
$$(f\,,g) \mapsto \int_{M}f g \ dM, \hspace{10pt}
f\,,g \in L_{2}(M),$$
where $dM$ is the volume element.

\begin{notation}
Let $A$ and $B$ be smooth manifolds and let
$\mathbb{F}=\mathbb{R}$ or $\mathbb{F}=\mathbb{C}$. We denote by
$C^{\ell,k}(A,B;\mathbb{F})$ the set of all functions $f: A \times B \rightarrow
\mathbb{F}$ whose mixed partial derivatives up to $\ell$ times
with respect to $A$ and up to $k$ times with respect to $B$ exist
and are continuous on $A\times B$.
\end{notation}

\begin{definition}
Let $M$ be any connected unbordered Riemannian manifold ($M$ need not to be complete in this
definition). A continuous function $H=H(x,y,t): M\times M \times (0,\infty)\rightarrow \mathbb{R}$
is called a {\it fundamental solution to the heat equation} on $M$
if it belongs to $C^{2,1}(M\times M,(0,\infty);\mathbb{R})$ and
satisfies the following conditions.
\begin{equation*}
\begin{split}
 &(1)\ \frac{\partial H}{\partial t}=-\Delta_{x}H;
\\ &(2)\ H(x,y,t)=H(y,x,t);
\\ &(3)\ \lim\limits_{t\to 0}\int_{M}H(x,y,t)f(y)\,dM(y)=f(x).
\end{split}
\end{equation*}
Here $\Delta_{x}$ is the Laplace--Beltrami operator with respect
to the first argument of $H$, and the convergence in (3) is
locally uniform in $x$ for every continuous function $f$ with
compact support on $M$.
\end{definition}

\begin{theorem}\label{uniqueness}
Let $M$ be any $m$-dimensional compact connected Riemannian
manifold without boundary. Then $M$ has a unique fundamental
solution $H_{M}$ of the heat equation and
\begin{equation*}
\begin{split}
&(1)\ H_{M}\in C^{\infty}(M\times M\times(0,\infty));\\
&(2)\ 0\leqslant
H_{M}(x,y,t)\leqslant c_{M}t^{-m/2}, \hspace{10pt} 0<t<1,
\end{split}
\end{equation*}
\noindent
where the constant $c_{M}$ depends on $M$.
\end{theorem}
\vspace{0.1cm}

\begin{definition}
The unique fundamental solution $H_{M}=H_{M}(x,y,t)$ of the heat
equation on $M$
is called the {\it heat kernel}
of $M$.
\end{definition}
\vspace{0.1cm}

\begin{theorem}[Spectral theorem]\label{spectral}
Let $M$ be a compact connected Riemannian manifold without
boundary. The eigenvalue problem $$\Delta\varphi=\lambda\varphi$$
\noindent has a complete orthonormal system of
$C^{\infty}$-eigenfunctions $\varphi_{0},\,\varphi_{1},\ldots$ in
$L_{2}(M)$ with corresponding eigenvalues
$\lambda_{0},\,\lambda_{1},\ldots$\, . These have the following
properties.
\begin{equation*}
\begin{split}
& (1)\ 0=\lambda_{0}<\lambda_{1}\leqslant \lambda_{2}\leqslant
\ldots
\, ,
\quad
\lambda_{n}\to \infty,\ n\to\infty;
\\ &(2)\
H_{M}(x,y,t)=\sum_{n=0}^{\infty}e^{-\lambda_{n}t}\varphi_{n}(x)\varphi_{n}(y),
\end{split}
\end{equation*}
\noindent where the series converges uniformly on $M\times M$ for
each $t>0$.
\end{theorem}
\vspace{0.2cm}

\begin{definition}[Spectrum of the manifold]
The set $\{\lambda_n\}_{n=0}^{\infty}$ of eigenvalues of the
Laplace--Beltrami operator $\Delta=-\text{div grad}$ on Riemannian
manifold $M$, together with the multipli\-cities with which they
occur, is called the {\it spectrum of manifold} $M$.
\end{definition}

\begin{definition}
Let $f\colon M\to \mathbb{R}$ be a continuous function. A continuous
function $u~=~u(x,t):M\times[0,\infty)\to \mathbb{R}$ is called a
\it{solution of the heat equation for the initial condition}
$u(x,0)~=~f(x),$\ \rm if $u\in C^{2,1}(M,(0,\infty);\mathbb{R})$ and
if $u$ satisfies the \it{heat equation}
\begin{equation*}
\begin{split}
&1.\ \frac{\partial u}{\partial t}=-\Delta_{x}u;\\ &2.\
u(x,0)=f(x),\hspace{10pt} x \in M.
\end{split}
\end{equation*}
\end{definition}

\begin{lemma}
Let $H_{M}=H_{M}(x,y,t)$ be a fundamental solution to the heat equation
on $M$ and let $f\colon M\to \mathbb{R}$ be a continuous function.
The function
$$u(x,t):=\int_{M}H_{M}(x,y,t)f(y)\,dM(y), \hspace{10pt} t>0,$$ has
a continuous extension to function
$u(x,t)\colon M\times[0,\infty)\to \mathbb{R}$,
which is a solution to the heat equation with initial condition $u(x,0)=f(x)$.
\end{lemma}


\section{The proof of Theorem A}\label{proofA}

Let $M_1,\ldots,M_6$ denote the 6 compact orientable flat
3-manifolds and $N_1,\ldots,N_4$ denote the non-orientable ones in
the order given in section \ref{prelim}. Fundamental sets of manifolds
$M_1,\ldots,M_6$, $N_1,\ldots,N_4$ are shown in section 6, Fig. 2 -- 11
correspondingly.

Suppose that $f:M'\to M$ is a covering. Such a manifold $M'$ is
determined by a subgroup $\pi_1(M')$ of the fundamental group
$\pi_1(M)$. When $\pi_1(M')$ is a normal subgroup of $\pi_1(M)$,
$f:M'\to M$ is said to be a {\it regular covering}, and the
quotient group $G~=~\pi_1(M)/\pi_1(M')$ acts on $M'$ as a group of
isometries of $M'$. The group $G$ is called the {\it covering
group}. In order to prove \mbox{Theorem A} we need the following

\begin{lemma}[\mbox{\rm \cite[p. 826]{Brooks}}]
\label{Brooks}
Let
$f:M'\to M$
be a regular covering
and
$G=\pi_1(M)/\pi_1(M')$
be a covering group.
Then
$$ H_{M}(x,y,t)
=
\sum_{g\in G}
H_{M'}(\Tilde{x},g \Tilde{y},t),$$
where
$\Tilde{x}$ and $\Tilde{y}$ are
any inverse images of
$x$
and $y$
under the covering $f:M'\to M$.
\end{lemma}

\begin{proof}
It is evident that the right-hand side satisfies all the
conditions of the definition of the fundamental solution to the
heat equation. The lemma now follows from the uniqueness statement
(Theorem \ref{uniqueness}).
\end{proof}

In the following,
$\Lambda=\{\boldsymbol{\ell}=m_1\boldsymbol{\ell}_{1}+m_2\boldsymbol{\ell}_{2}
+m_3\boldsymbol{\ell}_{3},\ \ m_i\in \mathbb{Z}\,,\
\boldsymbol{\ell}_{i}\in \mathbb{R}^3,\, i=1,2,3\}$ will always be
a lattice
for $3$-dimensional flat torus $M_1=\mathbb{R}^3/\Lambda$,
$\text{vol}(M)$ --- the volume of
the manifold $M$ with respect to the Riemannian metric and
$M^{\ast}$ will always denote the fundamental set of the manifold $M$.

\begin{theoremA}
The trace function of each compact flat $3$-manifolds
$M_1,\ldots,M_6,N_1,\ldots,N_4$ can be computed explicitly:
$\text{tr}(H_{M_i})=F_i$ for each manifold $M_i$, $i=1,\ldots,6$
and $\text{tr}(H_{N_i})=P_i$ for each manifold $N_i$,
$i=1,\ldots,4$, where the functions $F_i$ and $P_j$ are given in
section~$6$, Table~$1$.
\end{theoremA}

\begin{proof}

Now we give here a sketch of the proof.
Assume that the fundamental set of the ma\-ni\-fold~$M$ is formed
by such vectors that the fundamental set of the flat torus $M_1$
regular covering the manifold $M$ is formed by vectors
$\boldsymbol{\ell}_{1}$, $\boldsymbol{\ell}_{2}$
and $\boldsymbol{\ell}_{3}$. Using a
relationship between the heat kernels on manifold $M$ and on
regular covering of $M$ (Lemma \ref{Brooks}) and properties of the
Euler--Poisson integral $\int\limits_{-\infty}^{\infty}e^{-t^2}\,dt=\sqrt{\pi}$
we obtain an explicit formula for the trace function
$\text{tr}(H_{M})$.

We prove Theorem A in case of $M_1$ and $M_2$. Other cases are proved
similarly.

\subsection{Manifold $M_1$ (flat torus)}

\begin{pro}\label{pro m1}
The heat kernel of $M_1$, where $ y = x $, has the following expression.
$$
H_{M_1}(x,x,t)
=
\frac {1} {(4\pi t)^{\frac{3}{2}}}
\sum_{\bell \in \Lambda}
e^{-|\bell\,|^2/4t}
,
\hspace{10pt}
x \in M_1.
$$
\end{pro}

\begin{proof}
Let us consider the universal covering of the flat torus
$\mathbb{R}^{3}\stackrel{\infty}{\rightarrow}M_1$.
It is known that
$$
H_{\mathbb{R}^{3}}(x,y,t)
=
\frac{1}{(4\pi t)^{\frac{3}{2}}}
e^{-|x-y|^{2}/4t},
\hspace{10pt}
x,y
\in
\mathbb{R}^{3}.
$$

It follows from Lemma \ref{Brooks} that
$$
H_{M_1}(x, x, t)
=
\frac{1}{(4\pi t)^{\frac{3}{2}}}
\sum_{\gamma\in\pi_{1}(M_1)}
e^
{-|\widetilde{x} - \gamma\, \widetilde{x}|^{2}/4t}
=
\frac{1}{(4\pi t)^{\frac{3}{2}}}
\sum_{\bell \in\Lambda}
e
^{-|\bell|^{2}/4t},
$$
where $\widetilde{x} \in \mathbb{R}^3$ is any inverse image
of $x \in M_1$ under the covering
$\mathbb{R}^{3}\stackrel{\infty}{\rightarrow}M_1$.

Proposition \ref{pro m1} is now proved.
\end{proof}

Then using Proposition \ref{pro m1} we have
$$
\text{tr}(H_{M_1})=\int_{M_1} H_{M_1}(x,x,t)\,dx
=
\frac{\text{vol}(M_1)}{(4\pi t)^{\frac{3}{2}}}
\sum_{\bell \in\Lambda}
e
^{-|\bell|^{2}/4t}.
$$

This proves Theorem A in case of $M_1$.

\begin{rem}
Theorem A in case of $3$-dimensional flat torus $M_1$
is proved by Brooks \cite{Brooks}.
\end{rem}

\subsection{Manifold $M_2$}

Assume that the fundamental set of the manifold $M_2$ is formed by
linearly independent vectors $\boldsymbol{\ell}_{1}/2$,
$\boldsymbol{\ell}_{2}$ and $\boldsymbol{\ell}_{3}$ (see \S 6,
Fig. 3). From Theorem \ref{orientable} in case of $M_2$ we see
that $\boldsymbol{\ell}_1$ is orthogonal to $\boldsymbol{\ell}_2$
and $\boldsymbol{\ell}_3$. Let $\varphi$ be an angle between
$\bell_2$ and $\bell_3$, where $0 < \varphi < \pi $.

\begin{pro}\label{pro m2}
The formula for the heat kernel of $M_2$, where $y=x$ is as
follows.
$$
H_{M_2}(x,x,t)
=
\frac {1} {(4\pi t)^{\frac{3}{2}}}
\Big(
\sum_{\bell \in \Lambda}
e^{-|\bell\,|^2/4t}
+
\sum_{\bell \in \Lambda}
e^{-|z-\beta z+\bell\,|^2/4t}
\Big)
,
$$
where
$x \in M_2$
and
$z$
is any inverse image of $x$ under the coverings
$\mathbb{R}^3 \stackrel{\infty}{\rightarrow} M_1 \stackrel{2}{\rightarrow}M_2$
and
$\beta$ acts on $\mathbb{R}^3$ as follows
$$
\beta\, z=A z+\bell_1/2,
\hspace{10pt}
z \in \mathbb{R}^3,
$$
where
\begin{equation*}
A=
\begin{pmatrix}
1&0&0\\
0&-1&0\\
0&0&-1
\end{pmatrix}
, \hspace{10pt} \bell_1/2=
\begin{pmatrix}
\ell_1/2 \\ 0\\ 0
\end{pmatrix}.
\end{equation*}
Here $A$ is a matrix with respect to the orthonormal
basis
$\boldsymbol{e}_1=(1,0,0)$, $\boldsymbol{e}_2=(0,1,0)$,
$\boldsymbol{e}_3=(0,0,1)$ and $\bell_1/|\bell_1|=\boldsymbol{e}_1$.

\end{pro}

\begin{proof}
Now we consider the two-fold covering
$M_1\stackrel{2}{\rightarrow}M_2$
and covering group
$G=\pi_1(M_2)/\pi_1(M_1)$.
Observe that
$$\pi_{1}(M_2)=\pi_{1}(M_1)\cup\alpha\,\pi_{1}(M_1),$$
where $\alpha$ acts on $M_1$
as follows (see Theorem \ref{orientable} in case of $M_2$)
$$\alpha\, \tilde{x} = A \tilde{x} + \xi, \hspace{10pt}
\tilde{x} \in M_1,$$ where
\begin{equation*}
A=
\begin{pmatrix}
1&0&0\\ 0&-1&0\\ 0&0&-1
\end{pmatrix}
, \hspace{10pt} \xi=
\begin{pmatrix}
\ell_1/2 \\ 0\\ 0
\end{pmatrix}.
\end{equation*}
Here $A$ is a matrix with respect to the orthonormal basis
$\boldsymbol{e}_1=(1,0,0)$, $\boldsymbol{e}_2=(0,1,0)$,
$\boldsymbol{e}_3=(0,0,1)$.
Geometrically, $\alpha$ identifies
left and right faces of the fundamental set of $M_2$
with rotation through $\pi$ (Fig. 3).

From Lemma \ref{Brooks} we see
$$
H_{M_2}(x,x,t)
=
\sum_{g\in G}
H_{M_1}(\Tilde{x},g\Tilde{x},t)
=
H_{M_1}(\Tilde{x},\Tilde{x},t)
+
H_{M_1}(\Tilde{x},\alpha\Tilde{x},t),
$$
where
$\Tilde{x}$
is
any inverse image of $x$ under the covering
$M_1\stackrel{2}{\rightarrow}M_2$.

Let us consider the universal covering of the flat torus
$\mathbb{R}^{3}\stackrel{\infty}{\rightarrow}M_1$.
From Lemma \ref{Brooks} and Proposition \ref{pro m1} we have

$$
H_{M_1}(\widetilde{x},\widetilde{x},t)
=
\frac{1}{(4\pi t)^{\frac{3}{2}}}
\sum_{\bell \in\Lambda}
e
^{-|\bell|^{2}/4t},
$$
\begin{multline*}
H_{M_1}(\widetilde{x},\alpha\widetilde{x},t)
=
\frac{1}{(4\pi t)^{\frac{3}{2}}}
\sum_{\gamma\in\pi_{1}(M_1)}
e^
{-|\widetilde{(\widetilde{x})}-\gamma\, \widetilde{(\alpha\,\widetilde{x})}|^{2}/4t}\\
=
\frac{1}{(4\pi t)^{\frac{3}{2}}}
\sum_{\gamma\in\pi_{1}(M_1)}
e^
{-|\widetilde{(\widetilde{x})}-\widetilde{(\alpha\,
\widetilde{x})}+\widetilde{(\alpha\, \widetilde{x})}-\gamma\, \widetilde{(\alpha\, \widetilde{x})}|^{2}/4t}
=
\frac{1}{(4\pi t)^{\frac{3}{2}}}
\sum_{\bell \in\Lambda}
e
^{-|\widetilde{(\widetilde{x})} - \widetilde{ (\alpha\, \widetilde{x})}+\bell|^{2}/4t},
\end{multline*}
where
$ \widetilde{x}\in M_1$ and $\widetilde{(\widetilde{x})}\in \mathbb{R}^{3}$
is any inverse image of
$\widetilde{x}$
under the covering
$\mathbb{R}^{3}\stackrel{\infty}{\rightarrow}\mathbb{R}^{3}/\Lambda=M_1$.

Observe that
$
\widetilde{\alpha\widetilde{x}}
=
\widetilde{\alpha}\widetilde{\widetilde{x}}
$,
where $\widetilde{\alpha}$ acts on $\mathbb{R}^3$ as follows
$$
\widetilde{\alpha}\, \widetilde{\widetilde{x}}
=
A\widetilde{\widetilde{x}}+\xi,
\hspace{10pt}
\widetilde{\widetilde{x}} \in \mathbb{R}^3.
$$
Now letting
$z=\widetilde{\widetilde{x}}$ and $\beta = \widetilde{\alpha}$
proves Proposition \ref{pro m2}.
\end{proof}

From Proposition \ref{pro m2} we have
\begin{multline}\label{m2 one}
\text{tr}(H_{M_2})=\int_{M_2} H_{M_2}(x,x,t)\,dx\\
=\left(
\frac{1}{4\pi t} \right)^{\frac{3}{2}} \int_{M_2^{\ast}}
\sum_{\boldsymbol{\ell}\in\Lambda}e^{-|\boldsymbol{\ell}\,|^{2}/4t}dx + \left( \frac{1}{4\pi
t} \right)^{\frac{3}{2}}
\int_{M_2^{\ast}}\sum_{\boldsymbol{\ell}\in\Lambda}e^{-|x- \beta\, x+\boldsymbol{\ell}\,|^{2}/4t} dx.
\end{multline}

Throughout this proof we work with respect to the orthonormal
basis $\boldsymbol{e}_1$, $\boldsymbol{e}_2$ and $\boldsymbol{e}_3$.
Assume that $x=(x_1,x_2,x_3)$. The
following are easy to verify:

\begin{equation}\label{m2 two}
\begin{split}
&\boldsymbol{\ell}
=
(\ell_1 m,\ell_2 n+\ell_3 k \cos\varphi,k\ell_3\sin\varphi);
\ \
\beta
x
=
(x_1+\ell_{1}/2,-x_2,-x_3);
\\
&
|x-\beta x+\boldsymbol{\ell}\,|^2
=
(m+\frac{1}{2})^2\ell_1^2+(2 x_2+n\ell_2+k\ell_3\cos\varphi)^2
+
(2 x_3 + k\ell_3\sin\varphi)^2.
\end{split}
\end{equation}
Here we have used that
$\boldsymbol{\ell}_1/|\boldsymbol{\ell}_1|=\boldsymbol{e}_1$,
$\boldsymbol{\ell}_2/|\boldsymbol{\ell}_2|=\boldsymbol{e}_2$
and
$\varphi$ is an angle between $\bell_2$ and $\bell_3$.

Define the volume of $M_2$ as $\text{vol}(M_2)=\int_{M_2^{\ast}}dx$,
where $dx=dx_1dx_2dx_3$. Note that $\text{vol}(M_2)=\text{vol}(M_2^{\ast})$
since $M_2^{\ast}$ is a fundamental set of $M_2$.
It follows from (\ref{m2 one}) and
(\ref{m2 two})
that
\begin{multline*}
\text{tr}(H_{M_2})
=
\frac{\text{vol}(M_2)}{(4\pi t)^{3/2}}
\sum_
{\boldsymbol{\ell}\in\Lambda}
e^
{-|\boldsymbol{\ell}\,|^{2}/4t}
\\
+
\frac{\ell_1/2}{(4\pi t)^{3/2}}
\sum_
{m\in \mathbb{Z}}
e^{-(m+\frac{1}{2})^2\ell_1^2/4t}
\underbrace
{\sum_
{n,k\in \mathbb{Z}}
\iint
\limits_
{K}
e^{-[
(2 x_2+n\ell_2+k\ell_3\cos\varphi)^2
+
(2 x_3 + k\ell_3\sin\varphi)^2
]
/
4t}
dx_2 dx_3}_
{I_1}
,
\end{multline*}
where $K$ is a parallelogram spanned by vectors
$\boldsymbol{\ell}_2$ and $\boldsymbol{\ell}_3$.

Observe that integration in $I_1$ over $K$
can be changed by the integration over rec\-tangle~$R$
with sides of length $\ell_2$ and $\ell_3 \sin \varphi$.
Integral $I_1$
can be represented by Fubini's theorem as a product of integrals.
It is easy to verify using properties of the Euler-Poisson integral
$\int
\limits_
{-\infty}^
{\infty}
e^{-t^2}
dt
=
\sqrt{\pi}$
that

\begin{multline*}
I_1
=
\sum_
{k\in \mathbb{Z}}
\Biggl(
\int
\limits_
{0}^
{\ell_3 \sin\varphi}
e^
{-(2 x_3 + k \ell_3 \sin\varphi)^2/4t}
dx_3
\sum_
{n\in \mathbb{Z}}
\int
\limits_
{0}^
{\ell_2}
e^
{-(2 x_2 + n \ell_2 + k \ell_3 \cos\varphi)^2/4t}
dx_2
\Biggr)
\\
=
2 \sqrt{\pi t}
\sum_
{k\in \mathbb{Z}}
\int
\limits_
{0}^
{\ell_3 \sin\varphi}
e^
{-(2 x_3 + k \ell_3 \sin\varphi)^2/4t}
d x_3
=
4 \pi t.
\end{multline*}

Hence,
$$
\text{tr}
(H_{M_2})
=
\frac
{\text{vol}(M_2)}
{(4\pi t)^{3/2}}
\sum_
{\boldsymbol{\ell} \in \Lambda}
e^
{-|\boldsymbol{\ell}\,|^{2}/4t}
\\
+
\frac
{\ell_1}
{4 \sqrt{\pi t}}
\sum_
{m\in \mathbb{Z}}
e^
{-(m+\frac{1}{2})^2 \ell_1^2 / 4t}.
$$

Theorem A in case of $M_2$ is now proved.

\end{proof}


\section{The proof of Theorem B}\label{proofB}

\begin{theoremB}
Any two homeomorphic compact flat $3$-manifolds are isospectral if
and only if they are isometric.
\end{theoremB}\vspace{0.2cm}

The proof is based on the following propositions and Lemma
\ref{mainlemma}.
To describe the spectrum of the manifold $M$ we use the trace
function of $M$ $$\text{tr}(H_{M})=\int_{M}H_{M}(x,x,t)\,dM,$$
where $H_{M}(x,y,t)$ is a
fundamental solution to the heat equation on $M$, $dM$ is
the volume element.

\begin{pro}\label{pro1}
Let $M$ be a compact connected Riemannian manifold without
boundary. Then the trace
function $\text{tr}(H_{M})$ is determined exactly
by the spectrum of the ma\-ni\-fold $M$.
\end{pro}

\begin{proof}
By spectral theorem (Theorem \ref{spectral}) we have
$$\text{tr}(H_{M})=
\int\limits_{M}\sum_{n=0}^{\infty}
e^{-\lambda_{n}t}\varphi_{n}(x)\varphi_{n}(x)\,dM
=\sum_{n=0}^{\infty}
e^{-\lambda_{n}t}\int\limits_{M}\varphi_{n}(x)\varphi_{n}(x)\,dM
=\sum_{n=0}^{\infty}e^{-\lambda_{n}t},$$
since $\varphi_{0},\varphi_{1},\ldots$ are orthonormal in
$L_{2}(M)$. Here $\{\varphi_{n}\}_{n=0}^{\infty}$ is again a
complete orthonormal sequence of eigenfunctions of the
Laplace--Beltrami operator on the compact Riemannian ma\-ni\-fold $M$
without boundary and $0=\lambda_{0}<\lambda_{1}\leqslant \lambda_{2}\leqslant
\ldots$ is the corresponding sequence of eigenvalues.
Proposition \ref{pro1} is now proved.
\end{proof}

\begin{pro}\label{pro2}
Let $M$ be a compact connected Riemannian manifold without
boundary. Then the trace function $\text{tr}(H_{M})$ determines exactly
the spectrum of the manifold $M$.
\end{pro}

\begin{proof}
Let $\{\lambda_{n}\}_{n=0}^{\infty}$ be a sequence of eigenvalues
of the Laplace--Beltrami operator on $M$. By spectral theorem (Theorem \ref{spectral})
we have $0=\lambda_{0}<\lambda_{1}\leqslant \lambda_{2}\leqslant
\ldots$\ . Thus, $\lambda_{0}=0$ with \mbox{multiplicity $1$.} Assume that
$\lambda_{1},\ldots,\lambda_{k}$ have been found with their
multiplicities. Then from Proposition \ref{pro1} we see that $\lambda_{k+1}$
is the largest value $\lambda$ such that
$$F(\lambda)=\lim_{t\to +\infty}\frac{\text{tr}(H_{M})-
\sum\limits_{n=0}^{k}e^{-\lambda_{n}t}}{e^{-\lambda t}} <
\infty,$$
where $F(\lambda)$ is the multiplicity of $\lambda_{k+1}$.
This proves Proposition \ref{pro2}.
\end{proof}

\begin{rem}
Observe that the previous procedure in Proposition \ref{pro2}
of finding $\lambda_{k}$ in consecutive order is similar to the
procedure of finding coefficients in Taylor series. Here we find
coefficients at exponents and in contrast to Taylor series we find
powers of exponents. Also one can see that the proof of
Proposition \ref{pro2} is similar to the proof
of Huber's theorem, where a similar procedure is used
(see \cite[chapter 9]{Buser}).
\end{rem}\vspace{0.1cm}

\begin{lemma}\label{mainlemma}
Let $M$ be a compact flat $3$-manifold. Then the trace function
$\text{tr}(H_{M})$ determines up to isometry the corresponding
manifold $M$.
\end{lemma}

\begin{proof}
The proof is based on the following. Assume that the trace function of
the flat \mbox{$3$-manifold} $\text{tr}(H_{M})$ is given. The manifold $M$ is
determined up to isometry by its fundamental set. From Theorem \ref{isometry}
we know which parameters determine the fundamental
set of the corresponding manifold $M$. Then using a similar procedure as in
Proposition \ref{pro2} we find all these parameters from the trace
function $\text{tr}(H_{M})$.
Using Theorem~A we observe that
the trace function $\text{tr}(H_{M})$ determine a volume
of the corresponding manifold
$M$, namely,
$\text{Vol}(M) =
(4 \pi)^{\frac{3}{2}}\, \lim\limits_{ t \to +0 }
t^{ \frac{3}{2} }\, \text{tr}(H_{M})$.

Thus, from $\text{tr}(H_{M})$ we find the fundamental
set of the corresponding manifold $M$ and, consequently, we determine up to isometry
the manifold $M$.

\subsection{Manifold $\boldsymbol{M_2}$}

We assume as in section
\ref{proofA} that the fundamental set of the manifold $M_2$ is
formed by linearly independent vectors $\boldsymbol{\ell}_{1}/2$,
$\boldsymbol{\ell}_{2}$ and $\boldsymbol{\ell}_{3}$ (see \S 6, Fig. 3).
Observe that isometry classes of $M_2$ are
parameterized by $\text{vol}(M_2)$, $|\boldsymbol{\ell}_1|=\ell_1$
and $2$-dimensional lattice $ \Lambda_2 = \{\boldsymbol{\ell}_2,\,
\boldsymbol{\ell}_3\}$ generated by vectors $\boldsymbol{\ell}_2$ and
$\boldsymbol{\ell}_3$ (Theorem \ref{isometry}).
Now we find all these parameters
from the trace function $\text{tr}(H_{M_2})$.

(1) Using Theorem A we show that $\text{tr}(H_{M_2})$
determines exactly $\ell_1$. Let
$$
F(\lambda)= \lim_{t\to +0}
\frac{(4\pi t)^{3/2}\text{tr}(H_{M_2})-\text{vol}(M_2)}
{e^{-\lambda^{2}/4t}}.
$$
Hence,
$$
F(\lambda)
=
\lim_{ t \to +0 }
\Big(
\text{vol}(M_2)
\left(
2 e^{ ( \lambda^{2} - \ell_{1}^{2} ) / 4t }
+
2 e^{ ( \lambda^{2} - \ell_{2}^{2} ) / 4t }
+
\cdots
\right)
+ 2 \pi t\, \ell_{1}
\sum_{ m \in \mathbb{Z} }
e^{ ( \lambda^{2} - ( m + \frac{1}{2} )^{2} \ell_{1}^{2} ) / 4t }
\Big).
$$
We find $\lambda_{max}$ --- the largest positive value $\lambda$
such that $F(\lambda) < \infty$.
If \mbox{$F(\lambda_{max}) \neq 0$} then
$$
\text{tr}(H_{M_2}) := \text{tr}(H_{M_2}) -
2 \frac{\text{vol}(M_2)}{(4 \pi t)^{3/2}}
\sum_{k=1}^{\infty}
e^{- k^2 \lambda_{max}^{2} / 4t },
$$
again we find $\lambda_{max}$ --- the largest positive value $\lambda$
such that $F(\lambda) < \infty$.
If $\lambda_{max}$
is such that $F(\lambda_{max}) = 0$ then $\lambda_{max}=\ell_1 / 2$.

(2) Now we show that $\text{tr}(H_{M_2})$ determines
the plane lattice $\{\boldsymbol{\ell}_2,\,
\boldsymbol{\ell}_3\}$. Let
$$
F_1(\lambda)
=
\lim_{ t \to +0 }
\frac{(4\pi t)^{3/2}
\text{tr}(H_{M_2}) - \text{vol}(M_2) ( 1 + \sum\limits_{ m \in \mathbb{Z} }
e^{ - m^2 \ell_1^2 / 4t } )
-
2 \pi t \sum\limits_{ m \in \mathbb{Z} }
e^{ - ( m + \frac{1}{2} )^2 \ell_1^2 / 4t } }
{e^{ - \lambda^{2} / 4t } }.
$$

Hence,
\begin{multline*}
F_1( \lambda )
=
\lim_{ t \to +0 }
\Big(
\text{vol}(M_2)
\big(
2 e^{ ( \lambda^{2} - \ell_{2}^{2} ) / 4t }
+
2 e^{ ( \lambda^{2} - \ell_3^{2} ) / 4t }
+
4 e^{ ( \lambda^{2} - (\ell_1^2 + \ell_{2}^{2}) ) / 4t }
\\
+
4 e^{ ( \lambda^{2} - (\ell_1^2 + \ell_{3}^{2}) ) / 4t }
+
2 e^{ ( \lambda^{2} - (\ell_2^2 + \ell_{3}^{2} + 2 \ell_2 \ell_3 \cos\varphi) ) / 4t }
\\
+
2 e^{ ( \lambda^{2} - (\ell_2^2 + \ell_{3}^{2} - 2 \ell_2 \ell_3 \cos\varphi) ) / 4t }
+
4 e^{ ( \lambda^{2} - (\ell_1^2 + \ell_2^2 + \ell_{3}^{2} + 2 \ell_2 \ell_3 \cos\varphi) ) / 4t }
+
\cdots
\big)
\Big),
\end{multline*}
where $\varphi$ is an angle between $\boldsymbol{\ell}_2$ and $\boldsymbol{\ell}_3$,
$0 < \varphi < \pi $.

We find $\lambda_{max}$ --- the largest positive value $\lambda$
such that $F_1(\lambda_{max}) < \infty$ and then
$$
\text{tr}(H_{M_2})
:=
\text{tr}(H_{M_2})
-
2 \frac{\text{vol}(M_2)}{(4 \pi t)^{3/2}}
\sum_{k=1}^{\infty}
e^{- k^2 \lambda_{max}^{2} / 4t },
$$
again we find the next $\lambda_{max}$.
Using this procedure we determine that
$$\Lambda_2=\{\lambda\  |\  F_1(\lambda)=2\,\text{vol}(M_2)\}.$$

It is known \cite{Brooks} that the plane lattice
is determined up to
isometry by its shortest length, the shortest length of a linearly
independent lattice points, and the area of the parallelogram they
span. Since we can determine the area $S$ of the parallelogram spanned by
$\boldsymbol{\ell}_2$ and $\boldsymbol{\ell}_3$ as $S=\text{vol}(M_2)/(\ell_1/2)$
then the fundamental set of $M_2$ is determined
from the $\text{tr}(H_{M_2})$, this proving Lemma~\ref{mainlemma} in case of $M_2$.

\subsection{Manifold $\boldsymbol{M_4}$}

We assume that the fundamental set of the manifold $M_4$ is
formed by linearly independent vectors $\boldsymbol{\ell}_{1}/4$,
$\boldsymbol{\ell}_{2}$ and $\boldsymbol{\ell}_{3}$ (see \S
6, Fig. 5). Observe that isometry classes of $M_4$ are
parameterized by $\text{vol}(M_4)$ and
$\ell_1 = |\boldsymbol{\ell}_1|$ (Theorem
\ref{isometry}). Now we show that
$\text{tr}(H_{M_4})$ determines these parameters.

Using Theorem A we show that $\text{tr}(H_{M_4})$ determines
exactly $\ell_1$. Let $$F(\lambda)= \lim_{t\to +0}
\frac{(4\pi t)^{3/2}\text{tr}(H_{M_4})-\text{vol}(M_4)}
{e^{-\lambda^{2}/4t}}.$$
Hence,
\begin{multline*}
F( \lambda )
=
\lim_{ t \to +0 }
\Big(
\text{vol}(M_4)
\big(
2 e^{ ( \lambda^{2} - \ell_{1}^{2} ) / 4t }
+
4 e^{ ( \lambda^{2} - \ell^{2} ) / 4t }
+
\cdots
\big)
\\
+
2 \pi t \, \ell_{1}
\sum_{ m \in \mathbb{Z} }
e^{ ( \lambda^{2} - ( m + \frac{1}{4} )^{2} \ell_{1}^{2}) / 4t }
+
\pi t\, \ell_{1}
\sum_{ n \in \mathbb{Z} }
e^{ ( \lambda^{2} - ( n + \frac{1}{2} )^{2} \ell_{1}^{2} ) / 4t }
\Big).
\end{multline*}
We find $\lambda_{max}$ --- the largest positive value $\lambda$
such that $F(\lambda) < \infty$.

If
\mbox{ $ F( \lambda_{max} ) \neq 0 $ }
then
$$
\text{tr}(H_{M_4})
:=
\text{tr}(H_{M_4}) -
2 \frac{\text{vol}(M_4)}{(4 \pi t)^{3/2}}
\sum_{k=1}^{\infty}
e^{- k^2 \lambda_{max}^{2} / 4t },
$$
again we find $\lambda_{max}$ --- the largest positive value
$\lambda$ such that $F(\lambda) < \infty$.
If $F(\lambda_{max}) = 0$ then
we determine $\ell_1$ from
$\text{tr}(H_{M_4})$, moreover, observe that
$\lambda_{max} = \ell_1 /4$.

Thus, the trace function $\text{tr}(H_{M_4})$ determines
the fundamental set of $M_4$ and, therefore, the manifold $M_4$
up to isometry.
Lemma \ref{mainlemma} in case of $M_4$ is now proved.

\subsection{Manifold $\boldsymbol{M_6}$}

Assume that the fundamental set of the manifold $M_6$ is formed
by linearly independent vectors $\boldsymbol{\ell}_{1}/2$,
$\boldsymbol{\ell}_{2}/2$ and $\boldsymbol{\ell}_{3}$
(see \S 6, Fig. 7).
From Theorem \ref{orientable} in case of
$M_6$ we see that $\boldsymbol{\ell}_{i}, \, i = 1, 2, 3$ are mutually orthogonal
and let
$\boldsymbol{\ell}_{1} = (\ell_1, 0, 0)$,
\mbox{$\boldsymbol{\ell}_{2} = (0, \ell_2, 0)$,}
$\boldsymbol{\ell}_{3} = (0, 0, \ell_3)$
with respect to the orthonormal basis
$\boldsymbol{e}_1$, $\boldsymbol{e}_2$,
$\boldsymbol{e}_3$.

Observe that the fundamental set of 3-dimensional flat torus $M_1$
covering the manifold $M_6$ is formed by vectors
$\boldsymbol{\ell}_{1}$, $\boldsymbol{\ell}_{2}$ and
$\boldsymbol{\ell}_{3}$ (see \S 6, Fig. 7) and isometry classes of
$M_6$ are parameterized by the unordered set of three
$\{\ell_1,\,\ell_2,\,\ell_3\}$ \mbox{(Theorem \ref{isometry}).}
Now we find all these parameters from the trace function
$\text{tr}(H_{M_6})$.

Now we show using Theorem A that $\text{tr}(H_{M_6})$ determines
$\ell_1$, $\ell_2$ and $\ell_3$. Let $$ F(\lambda)= \lim_{t\to +0}
\frac{(4\pi t)^{3/2}\text{tr}(H_{M_6})-\text{vol}(M_6)}
{e^{-\lambda^{2}/4t}}.$$
Hence,
\begin{multline*}
F( \lambda )
=
\lim_{ t \to +0 }
\Big(
\text{vol}(M_6)
\big(
2 e^{ ( \lambda^{2} - \ell_{1}^{2} ) / 4t }
+
2 e^{ ( \lambda^{2} - \ell_{2}^{2} ) / 4t }
+
2 e^{ ( \lambda^{2} - \ell_{3}^{2} ) / 4t }
+
\cdots
\big)
\\
+
\pi t \,
\big(
\ell_{1}
\sum_{ m \in \mathbb{Z} }
e^{ ( \lambda^{2} - ( m + \frac{1}{2} )^{2} \ell_{1}^{2} ) / 4t }
\\
+
\ell_{2}
\sum_{ n \in \mathbb{Z} }
e^{ ( \lambda^{2} - ( n + \frac{1}{2} )^{2} \ell_{2}^{2} ) / 4t }
+
\ell_{3}
\sum_{ k \in \mathbb{Z} }
e^{ ( \lambda^{2} - ( k + \frac{1}{2} )^{2} \ell_{3}^{2} ) / 4t }
\big)
\Big).
\end{multline*}
We find $\lambda_{max}$ --- the largest positive value $\lambda$
such that $F(\lambda)<\infty$.

If $F(\lambda_{max}) \neq 0$ then
$$
\text{tr}(H_{M_6})
:=
\text{tr}(H_{M_6}) - 2 \frac{\text{vol}(M_6)}{(4 \pi t)^{3/2}}
\sum_{k=1}^{\infty}
e^{- k^2 \lambda_{max}^{2} / 4t },
$$
again we find $\lambda_{max}$.

If $F(\lambda_{max}) = 0$ then
$\lambda_{max} = \text{min}\{\ell_1 /2, \ell_2 /2, \ell_3 /2 \}$.
Recall that we find
$\ell_1,\,\ell_2$ and $\ell_3$ in any order (see Theorem \ref{isometry}).
Knowing $\lambda_{max}$, we put
$$
\text{tr}(H_{M_6})
:=
\text{tr}(H_{M_6}) - \frac{ 1 }{ 8\sqrt{\pi t} }
\sum\limits_{ m \in \mathbb{Z} }
e^{ - ( m + \frac{1}{2} )^{2} 4 \lambda_{max}^{2} / 4t }
$$
and then determine the second lowest
value of $\{\ell_1/2, \ell_2/2, \ell_3/2\}$
and so on. In this way we will find the unordered set of three
$\{\ell_1,\,\ell_2,\,\ell_3\}$.

Thus, the trace function $\text{tr}(H_{M_6})$ determines
the fundamental set of $M_6$ and, therefore, the manifold $M_6$
up to isometry.
Lemma \ref{mainlemma} in case of $M_6$ is now proved.

\subsection{Manifold $\boldsymbol{N_2}$}

Assume that the fundamental set
of the manifold $N_2$ is formed
by linearly independent vectors
$\boldsymbol{\ell}_{1}/2$,
$\boldsymbol{\ell}_{2}$ and
$\boldsymbol{\ell}_{3}$ (see \S
6, Fig. 9).
From Theorem \ref{nonorientable}
in case \mbox{of $N_2$} we see that
the orthogonal projection of
$\boldsymbol{\ell}_3$ on the
$(\boldsymbol{\ell}_1,\,\boldsymbol{\ell}_2)$ -- plane is
$(\boldsymbol{\ell}_{1}+\boldsymbol{\ell}_{2})/2$.
Let
\begin{equation*}
\begin{split}
&\boldsymbol{\ell}_{1}
=
(\ell_1,0,0),
\\
&\boldsymbol{\ell}_{2}
=
(\ell_2\cos\varphi,\ell_2\sin\varphi,0),
\\
&\boldsymbol{\ell}_{3}
=
(\frac{1}{2}(\ell_1+\ell_2\cos\varphi),\frac{1}{2}\ell_2\sin\varphi,h)
\end{split}
\end{equation*}
with respect to the orthonormal basis
$\boldsymbol{e}_1$, $\boldsymbol{e}_2$,
$\boldsymbol{e}_3$,
where
$\varphi$ is an angle between $\boldsymbol{\ell}_{1}$ and
$\boldsymbol{\ell}_{2}$,
$0 < \varphi < \pi $,
$h$ is a projection of $\boldsymbol{\ell_3}$ on
$\boldsymbol{e_3}$.

Observe that
the fundamental set of 3-dimensional flat torus $M_1$
covering the manifold $N_2$ is formed by vectors
$\boldsymbol{\ell}_{1}$, $\boldsymbol{\ell}_{2}$
and
$\boldsymbol{\ell}_{3}$ (see \S 6, Fig. 9)
and
isometry classes of $N_2$ are
parameterized by
$h$
and
the plane lattice
$\{\boldsymbol{\ell}_{1}, \boldsymbol{\ell}_{2}\}$
generated by vectors $\boldsymbol{\ell}_{1}$ and $\boldsymbol{\ell}_{2}$
(Theorem~\ref{isometry}).

Now we find as in previous cases all these parameters from the
trace function $\text{tr}(H_{N_2})$.
We will show that
$\text{tr}(H_{N_2})$
determines
$h=\text{vol}(N_2)/S$, where
$S
=
\ell_1\ell_2\sin\varphi$
and the plane lattice
$\{\boldsymbol{\ell}_{1}, \boldsymbol{\ell}_{2}\}$.
First,
we determine $S=\ell_1\ell_2\sin\varphi$ --- the area of the
parallelogram spanned by
$\boldsymbol{\ell}_{1}$ and
$\boldsymbol{\ell}_{2}$ and
then $h=\text{vol}(N_2)/S$.
Second,
we determine the plane lattice
$\{\boldsymbol{\ell}_{1}, \boldsymbol{\ell}_{2}\}$.

(1) Determine $S=\ell_1\ell_2\sin\varphi$ using the following limit

$$ F(\lambda)= \lim_{t\to +0}
\frac{(4\pi t)^{3/2}\text{tr}(H_{N_2})-\text{vol}(N_2)}
{e^{-\lambda^{2}/4t}}.$$
Hence,
\begin{multline*}
F( \lambda )
=
\lim_{ t \to +0 }
\Big(
\text{vol}(N_2)
\big(
2 e^{ ( \lambda^{2} - \ell_{1}^{2} ) / 4t }
+
2 e^{( \lambda^{2} - \ell_{2}^{2} ) / 4t }
+
2 e^{( \lambda^{2} - \ell_{3}^{2} ) / 4t}
+
\cdots
\big)
\\
+
\frac{ \sqrt{ \pi t } \ell_{1} \ell_{2} \sin\varphi }{ 2 }
\sum\limits_{ m,\,n \in \mathbb{Z} }
e^{   (\lambda^2 - (( \ell_1 ( m + \frac{1}{2} ) + n \ell_2 \cos\varphi)^2
+
( n \ell_2 \sin\varphi)^2)) / 4t  }
\\
+
\frac{ \sqrt{ \pi t } \ell_{1} \ell_{2} \sin\varphi }{ 2 }
\sum\limits_{ m,\,n \in \mathbb{Z}}
e^{   (\lambda^2 - ((\ell_1  m  + ( n + \frac{1}{2} ) \ell_2 \cos\varphi)^2
+
( ( n + \frac{1}{2} ) \ell_2 \sin\varphi)^2)) / 4t  }.
\end{multline*}
We find $\lambda_{max}$ --- the largest positive value $\lambda$
such that $F(\lambda)<\infty$.

If $F(\lambda_{max}) \neq 0$ then
$$
\text{tr}(H_{N_2})
:=
\text{tr}(H_{N_2}) - 2 \frac{\text{vol}(N_2)}{(4 \pi t)^{3/2}}
\sum_{ k = 1 }^{ \infty }
e^{ - k^2 \lambda_{max}^{2} / 4t },
$$
again we find $\lambda_{max}$.

If $F(\lambda_{max}) = 0$ then
fix this value $\lambda_{max}$ as $\lambda_{0}$.
Determine
$S=\ell_1\ell_2\sin\varphi$ using
the following limit
\begin{equation*}
2 \lim_{ t \to +0 }
\frac{ ( 4 \pi t )^{ 3 / 2 } \text{tr}(H_{N_2})
- \text{vol}(N_2) }
{ \sqrt{ \pi t }\, e^{ - \lambda_{0}^{2} / 4t } }
=
S,
\end{equation*}
and, therefore, $ h = 2 \text{vol}(N_2) / S $.

(2) Now we determine the plane lattice
$\{\boldsymbol{\ell}_{1}, \boldsymbol{\ell}_{2}\}$
generated by vectors $\boldsymbol{\ell}_{1}$ and $\boldsymbol{\ell}_{2}$.
Using Theorem \ref{nonorientable} in case of $N_2$ we observe
that $\boldsymbol{\ell}_3=\boldsymbol{h}
+\frac{1}{2}(\boldsymbol{\ell}_1+\boldsymbol{\ell}_2),$ where
$\boldsymbol{h}$ is orthogonal to $\boldsymbol{\ell}_1$ and
$\boldsymbol{\ell}_2$, $ \boldsymbol{h}=(0, 0, h) $.
Hence,
\begin{multline*}
S_1
=
\sum\limits_{ \boldsymbol{\ell} \in \Lambda }
e^{ - | \boldsymbol{\ell} |^2 / 4t }
=
\sum\limits_{ m, k, n \in \mathbb{Z} }
e^{ - | k \boldsymbol{\ell}_1 + m \boldsymbol{\ell}_2 + n \boldsymbol{\ell}_3 |^2 / 4t }
\\
=
\sum\limits_{ m, k, n \in \mathbb{Z} }
e^{ - | k \boldsymbol{\ell}_1 + m \boldsymbol{\ell}_2
+
\frac{n}{2}
( \boldsymbol{\ell}_1 + \boldsymbol{\ell}_2)
+
n \boldsymbol{h} |^2 / 4t }
=
\sum\limits_{ m, k, n \in \mathbb{Z} }
e^{ - ( | k \boldsymbol{\ell}_1 + m \boldsymbol{\ell}_2
+
\frac{n}{2}
( \boldsymbol{\ell}_1 + \boldsymbol{\ell}_2 ) |^2
+
n^2 h^2 ) / 4t }.
\end{multline*}

Now we split $S_1$ in two sums $S_{11}$ and $S_{12}$ corresponding to $n=2j$ and
$n=2j+1$.

\begin{multline*}
S_{11}
=
\sum\limits_{ m, k, n \in \mathbb{Z} }
e^{ - ( | k \boldsymbol{\ell}_1 + m \boldsymbol{\ell}_2
+
\frac{n}{2}
( \boldsymbol{\ell}_1 + \boldsymbol{\ell}_2 ) |^2 + n^2 h^2 ) / 4t }|_{ n = 2j }
=
\sum\limits_{ m, k, j \in \mathbb{Z} }
e^{ - ( | k \boldsymbol{\ell}_1 + m \boldsymbol{\ell}_2
+
\frac{2j}{2} ( \boldsymbol{\ell}_1 + \boldsymbol{\ell}_2 ) |^2 + ( 2j )^2 h^2 ) / 4t }
\\
=
\sum\limits_{ m, k, j \in \mathbb{Z} }
e^{ - ( | ( k + j ) \boldsymbol{\ell}_1 + ( m + j ) \boldsymbol{\ell}_2 |^2
+
( 2j )^2 h^2 ) / 4t }
=
\sum\limits_{ j \in \mathbb{Z} }
e^{ - j^2 h^2 / t }
\Big(
\sum\limits_{ m, k \in \mathbb{Z} }
e^{ - | ( k + j ) \boldsymbol{\ell}_1 + ( m + j ) \boldsymbol{\ell}_2 |^2 / 4t }
\Big)
\\
=
\sum\limits_{ j \in \mathbb{Z} }
e^{ - j^2 h^2 / t }
\Big(
\sum\limits_{ m, k \in \mathbb{Z} }
e^{ - | k \boldsymbol{\ell}_1 + m \boldsymbol{\ell}_2 |^2 / 4t }
\Big)
=
\sum\limits_{ m, k \in \mathbb{Z} }
e^{ - | k \boldsymbol{\ell}_1 + m \boldsymbol{\ell}_2 |^2 / 4t }
+
2 \sum\limits_{ j = 1 }^{ + \infty }
e^{ - j^2 h^2 / t }
\Big(
\sum\limits_{ m, k \in \mathbb{Z}}
e^{ - | k \boldsymbol{\ell}_1 + m \boldsymbol{\ell}_2 |^2 / 4t }
\Big).
\end{multline*}

\begin{multline*}
S_{12}
=
\sum\limits_{ m, k, n \in \mathbb{Z} }
e^{ - ( | k \boldsymbol{\ell}_1 + m \boldsymbol{\ell}_2
+
\frac{n}{2}
( \boldsymbol{\ell}_1 + \boldsymbol{\ell}_2 ) |^2 + n^2 h^2 ) / 4t }|_{ n = 2j+1 }
=
\sum\limits_{ m, k, j \in \mathbb{Z} }
e^{ - ( | k \boldsymbol{\ell}_1 + m \boldsymbol{\ell}_2
+
\frac{ 2j + 1 }{2}
( \boldsymbol{\ell}_1 + \boldsymbol{\ell}_2 ) |^2 + ( 2j + 1 )^2  h^2 ) / 4t }
\\
=
\sum\limits_{ m, k, j \in \mathbb{Z} }
e^{ - ( | ( k + j ) \boldsymbol{\ell}_1 + ( m + j ) \boldsymbol{\ell}_2
+
\frac{1}{2}
( \boldsymbol{\ell}_1 + \boldsymbol{\ell}_2 ) |^2 + ( 2j + 1 )^2 h^2 ) / 4t }
\\
=
\sum\limits_{ j \in \mathbb{Z} }
e^{ - ( 2j + 1 )^2 h^2 / 4t }
\Big(
\sum\limits_{ m, k \in \mathbb{Z} }
e^{ - | ( k + j ) \boldsymbol{\ell}_1 + ( m + j ) \boldsymbol{\ell}_2
+
\frac{1}{2}
( \boldsymbol{\ell}_1 + \boldsymbol{\ell}_2 ) |^2 / 4t }
\Big)
\\
=
\sum\limits_{ j \in \mathbb{Z} }
e^{ - ( 2j + 1 )^2 h^2 / t }
\Big(
\sum\limits_{ m, k \in \mathbb{Z} }
e^{ - | k \boldsymbol{\ell}_1 + m \boldsymbol{\ell}_2
+
\frac{1}{2}
( \boldsymbol{\ell}_1 + \boldsymbol{\ell}_2 ) |^2 / 4t }
\Big)
\\
=
2 \sum\limits_{ j = 0 }^{ + \infty }
e^{ - ( 2j + 1 )^2 h^2 / t }
\Big(
\sum\limits_{ m, k \in \mathbb{Z} }
e^{ - | k \boldsymbol{\ell}_1 + m \boldsymbol{\ell}_2
+
\frac{1}{2}
( \boldsymbol{\ell}_1 + \boldsymbol{\ell}_2 ) |^2 / 4t }
\Big).
\end{multline*}

Thus, we have the following

\begin{multline*}
\text{tr}(H_{N_2})
=
\frac{\text{vol}(N_2)}{(4 \pi t)^{3/2}}
\sum\limits_{ m, k \in \mathbb{Z} }
e^{ - | k \boldsymbol{\ell}_1 + m \boldsymbol{\ell}_2 |^2 / 4t }
+
\frac{ 2 \text{vol}(N_2) }{ ( 4 \pi t )^{3/2} }
\sum\limits_{ j = 1 }^{ + \infty }
e^{ - j^2 h^2 / t }
\Big(
\sum\limits_{ m, k \in \mathbb{Z} }
e^{ - | k \boldsymbol{\ell}_1 + m \boldsymbol{\ell}_2 |^2 / 4t }
\Big)
\\
+
\frac{ 2 \text{vol}(N_2) }{ ( 4 \pi t )^{3/2} }
\sum\limits_{ j = 0 }^{ + \infty }
e^{ - ( 2j + 1 )^2 h^2 / t }
\Big(
\sum\limits_{ m, k \in \mathbb{Z} }
e^{ - | k \boldsymbol{\ell}_1 + m \boldsymbol{\ell}_2
+
\frac{1}{2}
( \boldsymbol{\ell}_1 + \boldsymbol{\ell}_2 ) |^2 / 4t }
\Big)
\\
+
\frac{ \ell_{1} \ell_{2} \sin\varphi }{ 8 \sqrt{\pi t} }
\sum\limits_{ m,\,n \in \mathbb{Z} }
e^{ - ( ( \ell_1 ( m + \frac{1}{2} ) + n \ell_2 \cos\varphi)^2
+
( n \ell_2 \sin\varphi)^2 ) / 4t }
\\
+
\frac{ \ell_{1} \ell_{2} \sin\varphi }{ 8 \sqrt{\pi t} }
\sum\limits_{ m,\,n \in \mathbb{Z} }
e^{ - ( ( \ell_1  m  + ( n + \frac{1}{2} ) \ell_2 \cos\varphi)^2
+
( ( n + \frac{1}{2} ) \ell_2 \sin\varphi)^2 ) / 4t  }
\end{multline*}

In order to determine the plane lattice
$\{\boldsymbol{\ell}_{1}, \boldsymbol{\ell}_{2}\}$
we consider the limit
$$ F(\lambda)= \lim_{t\to +0}
\frac{(4\pi t)^{3/2}\text{tr}^{\ast}(H_{N_2})}
{e^{-\lambda^{2}/4t}},$$
where
$$
\text{tr}^{\ast}(H_{N_2})=\text{tr}(H_{N_2})-
\frac{\text{vol}(N_2)}{(4 \pi t)^{3/2}}
-\frac{2\text{vol}(N_2)}{(4 \pi t)^{3/2}}
\sum\limits_{j=1}^{+\infty}e^{-j^2h^2/t}.
$$
We find $\lambda_{max}$ --- the largest positive value $\lambda$
such that $F(\lambda)<\infty$.

Consider the case $ F( \lambda_{max} ) \neq 0 $.
If
$$
F(\lambda_{max})
=
\frac{2 \text{vol}(N_2)}{(4 \pi t)^{3/2}}
\sum_{ j = 0 }^{ \infty }
e^{ - ( 2j + 1 )^2 h^2 / t}
$$
then
$\lambda_{max}
=\min\limits_{ m,\, n \in \mathbb{Z} }
\{ m \boldsymbol{\ell}_1 + n \boldsymbol{\ell}_2
+
\frac{1}{2}( \boldsymbol{\ell}_1 + \boldsymbol{\ell}_2 ) \}$,
otherwise,
$
\lambda_{max}
=
\min\limits_{ m,\, n \in \mathbb{Z} }
\{ m \boldsymbol{\ell}_1 + n \boldsymbol{\ell}_2 \} $.
Knowing $\lambda_{max}$, we put
in the first case
$$
\text{tr}^{\ast}(H_{N_2})
:=
\text{tr}^{\ast}(H_{N_2}) - \frac{ 2 \text{vol}(N_2) }{ (4 \pi t)^{3/2} }
\sum\limits_{ j = 0 }^{ \infty }
e^{ - ( 2j + 1 )^2 h^2 / t }
e^{ - \lambda_{max}^{2} / 4t },
$$
and in the latter case
$$
\text{tr}^{\ast}(H_{N_2})
:=
\text{tr}^{\ast}(H_{N_2}) - \frac{ \text{vol}(N_2)  }{ (4 \pi t)^{3/2} }
e^{ - \lambda_{max}^{2} / 4t }
+
\frac{ 2 \text{vol}(N_2) }{ (4 \pi t)^{3/2} }
\sum\limits_{ j = 1 }^{ \infty }
e^{ - j^2 h^2 / t }
e^{ - \lambda_{max}^{2} / 4t }.
$$
Then we determine the next value
$\lambda_{max}$.

If
$ F(\lambda_{max}) = 0 $
then put
$$
\text{tr}^{\ast}(H_{N_2})
:= \text{tr}^{\ast}(H_{N_2}) - \frac{\ell_1 \ell_2 \sin \varphi }{ 8 \sqrt{\pi t} }
\sum\limits_{ m, n, k \in \mathbb{Z} }
e^{ - \lambda_{max}^2 / 4t}
$$
and determine the next value
$\lambda_{max}$.

Thus, we obtain the plane lattice
\mbox{
$\{\boldsymbol{\ell}_1,\,\boldsymbol{\ell}_2\}
=
\{ m \boldsymbol{\ell}_1 + n \boldsymbol{\ell}_2,\, m, n \in \mathbb{Z} \}$.
}
To complete the proof in case of $N_2$ we observe that
a lattice
$\{\boldsymbol{\ell}_1,\,\boldsymbol{\ell}_2\}$ in $\mathbb{R}^2$
is determined up to
isometry by its shortest length, the shortest length of a linearly
independent lattice points, and the area of the parallelogram they
span. Since all of these parameters are determined from $\text{tr}(H_{N_2})$,
this proving Lemma \ref{mainlemma} in case of $N_2$.

\subsection{Manifold $\boldsymbol{N_3}$}
Assume that the fundamental set of the manifold $N_3$ is
formed by linearly independent vectors $\boldsymbol{\ell}_{1}/2$,
$\boldsymbol{\ell}_{2}/2$ and $\boldsymbol{\ell}_{3}$ (see \S
6, Fig. 10).
From Theorem \ref{nonorientable}
in case \mbox{of $N_3$} we see that $\boldsymbol{\ell}_{1}$,
$\boldsymbol{\ell}_{2}$ and $\boldsymbol{\ell}_{3}$ are
mutually orthogonal. Let $\boldsymbol{\ell}_{1}=(\ell_1,0,0)$,
$\boldsymbol{\ell}_{2}=(0,\ell_2,0)$
and $\boldsymbol{\ell}_{3}=(0,0,\ell_3)$
with respect to the orthonormal
basis $\boldsymbol{e}_1$, $\boldsymbol{e}_2$,
$\boldsymbol{e}_3$.
Observe that isometry classes of $N_3$ are
parameterized by $\text{vol}(N_3)$, $\ell_1$ and
$\ell_2$ or $\ell_3$.
(Theorem \ref{isometry}).

Now we find as in previous cases all these parameters from the
trace function $\text{tr}(H_{N_3})$, namely, we will show that
$\text{tr}(H_{N_3})$ determines
$\ell_1$, $\ell_2$ and $\ell_3$.

Let
$$
F(\lambda)
=
\lim_{ t \to +0 }
\frac{ (4 \pi t)^{3/2}\, \text{tr}(H_{N_3}) - \text{vol}(N_3) }
{ e^{ - \lambda^{2} / 4t } }.
$$
Hence,
\begin{multline*}
F(\lambda)
=
\lim_{ t \to +0 }
\Big(
\text{vol}(N_3)
\left(
2 e^{ ( \lambda^{2} - \ell_{1}^{2} ) / 4t }
+
2 e^{ ( \lambda^{2} - \ell_2^{2} ) / 4t }
+
2 e^{ ( \lambda^{2} - \ell_3^{2} ) / 4t }
+
\cdots
\right)
\\
+
\ell_{1} \pi t
\sum_{ m \in \mathbb{Z} }
e^{ ( \lambda^{2} - ( m + \frac{1}{2} )^{2} \ell_{1}^{2} ) / 4t }
+
\frac{ \ell_{1} \ell_2 \sqrt{\pi t} }{2}
\sum_{ m, n \in \mathbb{Z} }
e^{ ( \lambda^{2} - ( m^2 \ell_1^2 + ( n + \frac{1}{2} )^{2} \ell_{2}^{2} ) ) / 4t }
\\
+
\frac{ \ell_{1} \ell_3 \sqrt{\pi t} }{2}
\sum_{ m, k \in \mathbb{Z} }
e^{ ( \lambda^{2} - ( ( m + \frac{1}{2} )^{2} \ell_{1}^{2} + k^2 \ell_3^2 ) ) / 4t }
\Big).
\end{multline*}
Find $\lambda_{max}$ --- the largest positive value $\lambda$
such that $F(\lambda)<\infty$.
We have the following cases: (i) $F(\lambda_{max}) \neq 0 $ and
(ii) $F(\lambda_{max}) = 0$.

Consider case (i). In this case $\lambda_{max} = \ell_3$. Consequently, we have
\mbox{$\ell_1 \ell_2 = 4 \text{vol}(N_3) / \ell_3$.}
Let
$$
F_1(\lambda) =
\lim_{ t \to +0 }
\frac{ \pi t\, \text{tr}(H_{N_3}) - \frac{\text{vol}(N_3)}{8 \sqrt{ \pi t}}
\sum\limits_{ m \in \mathbb{Z} }
e^{ - m^2 \ell_3^2 / 4t} }
{ e^{ - \lambda^{2} / 4t } }.
$$
Hence,
\begin{multline*}
F_1(\lambda)
=
\lim_{ t \to +0 }
\Big(
\frac{ \text{vol}(N_3) } {8 \sqrt{\pi t } }
\left(
2 e^{ ( \lambda^{2} - \ell_{1}^{2} ) / 4t }
+
2 e^{ ( \lambda^{2} - \ell_2^{2} ) / 4t }
+
4 e^{ ( \lambda^{2} - ( \ell_1^2 + \ell_3^{2} ) ) / 4t }
+
\cdots
\right)
\\
+
\frac{ \ell_{1} \sqrt{\pi t} }{8}
\sum_{ m \in \mathbb{Z} }
e^{ ( \lambda^{2} - ( m + \frac{1}{2})^{2} \ell_{1}^{2} ) / 4t }
+
\frac{ \ell_{1} \ell_2 }{16}
\sum_{ m, n \in \mathbb{Z} }
e^{ ( \lambda^{2} - ( m^2 \ell_1^2 + ( n + \frac{1}{2} )^{2} \ell_{2}^{2} ) ) / 4t }
\\
+
\frac{ \ell_{1} \ell_3 }{16}
\sum_{ m, k \in \mathbb{Z} }
e^{ ( \lambda^{2} - ( ( m + \frac{1}{2} )^{2} \ell_{1}^{2} + k^2 \ell_3^2 ) ) / 4t }
\Big).
\end{multline*}
We find $\lambda_{max}$ --- the largest positive value
$\lambda$ such that $F_1(\lambda)<\infty$.

Observe that
$\lambda_{max} = \min \{ \ell_1 / 2 , \ell_2 / 2\}$. If
$F_1(\lambda_{max}) = \frac{ \ell_1 \ell_2}{16}$ then
$\lambda_{max} = \ell_2 / 2$, otherwise,
\mbox{$\lambda_{max} = \ell_1 / 2$.}
Thus, in case~(i) the trace function $\text{tr}(H_{N_3})$
determines all parameters of the fundamental set of
$N_3$.

Consider case (ii). In this case
$\lambda_{max} = \min \{ \ell_1 / 2, \ell_2 / 2 \}$. Let
$$
F(\lambda) =
\lim_{ t \to +0 }
\frac{ \pi t\, \text{tr}(H_{N_3}) - \frac{\text{vol}(N_3)}{ 8 \sqrt{ \pi t} } }
{ e^{ - \lambda^{2} / 4t } }.
$$
Hence,
\begin{multline*}
F(\lambda)
=
\lim_{ t \to +0 }
\Big(
\frac{ \text{vol}(N_3)} {8 \sqrt{\pi t }}
\left(
2 e^{(\lambda^{2} - \ell_{1}^{2}) / 4t}
+
2 e^{(\lambda^{2} - \ell_2^{2}) / 4t}
+
2 e^{(\lambda^{2}- \ell_3^{2}) / 4t}
+
\cdots
\right)
\\
+
\frac{ \ell_{1} \sqrt{\pi t} }{8}
\sum_{ m \in \mathbb{Z} }
e^{ ( \lambda^{2} - ( m + \frac{1}{2} )^{2} \ell_{1}^{2} ) / 4t }
+
\frac{ \ell_{1} \ell_2 }{16}
\sum_{ m, n \in \mathbb{Z} }
e^{ ( \lambda^{2} - ( m^2 \ell_1^2 + ( n + \frac{1}{2} )^{2} \ell_{2}^{2} ) ) / 4t }
\\
+
\frac{ \ell_{1} \ell_3 }{16}
\sum_{ m, k \in \mathbb{Z} }
e^{ ( \lambda^{2} - ( ( m + \frac{1}{2} )^{2} \ell_{1}^{2} + k^2 \ell_3^2 ) ) / 4t }
\Big).
\end{multline*}
We find $\lambda_{max}$ --- the largest positive value
$\lambda$ such that $F(\lambda)<\infty$. Fix
$\lambda_{max}$ as $\lambda_{0}$ and $F(\lambda_{max})$
as $f_{0}$. We have two cases: (1) $\lambda_{0} = \ell_1 / 2$ and
$f_{0} = \frac{\ell_1 \ell_3}{ 16 }$;
(2)~$\lambda_{0} = \ell_2 / 2$ and
$f_{0} = \frac{\ell_1 \ell_2}{ 16 }$.
Let
$$
F_1(\lambda)
=
\lim_{ t \to +0 }
\frac{ (4 \pi t)^{3/2}\, \text{tr}(H_{N_3}) - \text{vol}(N_3) - A }
{ e^{ - \lambda^{2} / 4t } },
$$
where
$$
A
=
f_{0} \sum_{ m, k \in \mathbb{Z} }
e^{ - ( m^2 B^2 + ( 2 \lambda )^2 ( k + \frac{1}{2} )^2 ) / 4t },
\hspace{10pt}
B = \frac{16 f_{0}}{2 \lambda}.
$$

Consequently, in case (1) we have
\begin{multline*}
F_1(\lambda)
=
\lim_{ t \to +0 }
\Big(
\text{vol}(N_3)
\left(
2 e^{ ( \lambda^{2} - \ell_{1}^{2} ) / 4t }
+
2 e^{ ( \lambda^{2} - \ell_2^{2} ) / 4t }
+
2 e^{ ( \lambda^{2} - \ell_3^{2} ) / 4t }
+
\cdots
\right)
\\
+
\ell_{1} \pi t
\sum_{ m \in \mathbb{Z} }
e^{ ( \lambda^{2} - ( m + \frac{1}{2} )^{2}\ell_{1}^{2} ) / 4t }
+
\frac{ \ell_{1} \ell_2 \sqrt{\pi t} }{2}
\sum_{ m, n \in \mathbb{Z} }
e^{ ( \lambda^{2} - ( m^2 \ell_1^2 + ( n + \frac{1}{2} )^{2} \ell_{2}^{2} ) ) / 4t }
\Big).
\end{multline*}
In the latter case we have
\begin{multline*}
F_1(\lambda)
=
\lim_{ t \to +0 }
\Big(
\text{vol}(N_3)
\left(
2 e^{ ( \lambda^{2} - \ell_{1}^{2} ) / 4t }
+
2 e^{ ( \lambda^{2} - \ell_2^{2} ) / 4t }
+
2 e^{ ( \lambda^{2} - \ell_3^{2} ) / 4t }
+
\cdots
\right)
\\
+
\ell_{1} \pi t
\sum_{ m \in \mathbb{Z} }
e^{ ( \lambda^{2} - (m + \frac{1}{2})^{2} \ell_{1}^{2} ) / 4t }
+
\frac{ \ell_{1} \ell_3 \sqrt{\pi t} }{2}
\sum_{ m, k \in \mathbb{Z} }
e^{ ( \lambda^{2} - ( ( m + \frac{1}{2} )^{2} \ell_{1}^{2} + k^2 \ell_3^2 ) ) / 4t }
\Big).
\end{multline*}
Find $\lambda_{max}$ --- the largest positive value
$\lambda$ such that $F_1(\lambda)<\infty$.
If $F_1(\lambda_{max}) \neq 0$ then
$$
\text{tr}(H_{N_3})
:=
\text{tr}(H_{N_3}) - \frac{ 2 \text{vol}(N_3) }{(4 \pi t)^{3/2}}
\sum_{ k = 1 }^{ \infty }
e^{- k^2 \lambda_{max}^{2} / 4t },
$$
again we find $\lambda_{max}$.
If $F_1(\lambda_{max}) = 0$ then consider
$$
F_2(\lambda) = \lim_{t \to +0}
\frac{\pi t\,
\text{tr}(H_{N_3}) - \frac{\text{vol}(N_3)}{8 \sqrt{\pi t}}}
{e^{-\lambda^{2}/4t}}.
$$
Hence, in case (1) we have
\begin{multline*}
F_2(\lambda)
= \lim_{ t \to +0 }
\Big(
\frac{ \text{vol}(N_3) } {8 \sqrt{\pi t }}
\left(
2 e^{ ( \lambda^{2} - \ell_{1}^{2} ) / 4t }
+
2 e^{ ( \lambda^{2} - \ell_2^{2} ) / 4t }
+
2 e^{ ( \lambda^{2} - \ell_3^{2} ) / 4t }
+
\cdots
\right)
\\
+
\frac{ \ell_{1} \sqrt{\pi t} }{8}
\sum_{ m \in \mathbb{Z} }
e^{ ( \lambda^{2} - (m + \frac{1}{2})^{2} \ell_{1}^{2} ) / 4t }
+
\frac{ \ell_{1} \ell_2 }{16}
\sum_{ m, n \in \mathbb{Z} }
e^{ ( \lambda^{2} - ( m^2 \ell_1^2 + ( n + \frac{1}{2} )^{2} \ell_{2}^{2} ) ) / 4t }
\Big),
\end{multline*}
in case (2)
\begin{multline*}
F_2(\lambda)
=
\lim_{ t \to +0 }
\Big(
\frac{ \text{vol}(N_3) } {8 \sqrt{\pi t }}
\left(
2 e^{ ( \lambda^{2} - \ell_{1}^{2} ) / 4t }
+
2 e^{ ( \lambda^{2} - \ell_2^{2} ) / 4t }
+
2 e^{ ( \lambda^{2} - \ell_3^{2} ) / 4t }
+
\cdots
\right)
\\
+
\frac{ \ell_{1} \sqrt{\pi t} }{8}
\sum_{ m \in \mathbb{Z} }
e^{ ( \lambda^{2} - (m + \frac{1}{2})^{2} \ell_{1}^{2} ) / 4t }
+
\frac{ \ell_{1} \ell_3 }{16}
\sum_{ m, k \in \mathbb{Z} }
e^{ ( \lambda^{2} - ( ( m + \frac{1}{2} )^{2} \ell_{1}^{2} + k^2 \ell_3^2 ) ) / 4t }
\Big).
\end{multline*}
Find $\lambda_{max}$ --- the largest positive value
$\lambda$ such that $F_{2}(\lambda)<\infty$. Observe that
$\lambda_{max} = \ell_1 / 2$. Here we have used that
in the first case $\ell_1 = \min \{ \ell_1, \ell_2 \}$,
in the latter case $\ell_2 = \min \{ \ell_1, \ell_2 \}$.

Again we consider
$$
F(\lambda)
=
\lim_{ t \to +0 }
\frac{ \pi t\, \text{tr}(H_{N_3}) - \frac{\text{vol}(N_3)}{ 8 \sqrt{ \pi t} } }
{ e^{ - \lambda^{2} / 4t } }.
$$
Hence,
\begin{multline*}
F(\lambda)
=
\lim_{ t \to +0 }
\Big(
\frac{ \text{vol}(N_3) } {8 \sqrt{\pi t }}
\left(
2 e^{ ( \lambda^{2} - \ell_{1}^{2} ) / 4t }
+
2 e^{ ( \lambda^{2} - \ell_2^{2} ) / 4t }
+
2 e^{ ( \lambda^{2} - \ell_3^{2} ) / 4t }
+
\cdots
\right)
\\
+
\frac{ \ell_{1} \sqrt{\pi t} }{8}
\sum_{ m \in \mathbb{Z} }
e^{ ( \lambda^{2} - (m + \frac{1}{2})^{2} \ell_{1}^{2} ) / 4t }
+
\frac{ \ell_{1} \ell_2 }{16}
\sum_{ m, n \in \mathbb{Z} }
e^{ ( \lambda^{2} - ( m^2 \ell_1^2 + ( n + \frac{1}{2} )^{2} \ell_{2}^{2} ) ) / 4t }
\\
+
\frac{ \ell_{1} \ell_3 }{16}
\sum_{ m, k \in \mathbb{Z} }
e^{ ( \lambda^{2} - ( ( m + \frac{1}{2} )^{2} \ell_{1}^{2} + k^2 \ell_3^2 ) ) / 4t }
\Big).
\end{multline*}
Find $\lambda_{max}$ --- the largest positive value
$\lambda$ such that $F(\lambda)<\infty$.

Observe that
$\lambda_{max} = \min \{ \ell_1 / 2,\,\ell_2 / 2 \}$.
Since we know the value
$\ell_1$ then
$\ell_3 = \frac{16 F(\lambda_{max})}{\ell_1}$ if $\lambda_{max} = \ell_1 / 2$,
and
$\lambda_{max} = \ell_2 / 2$ if $\lambda_{max} \neq \ell_1 / 2$.

Thus, in case (ii) the trace function $\text{tr}(H_{N_3})$
determines all parameters of the fundamental set of
$N_3$.
Lemma~\ref{mainlemma} in case of~$N_3$ is proved.

\begin{rem}
Lemma \ref{mainlemma} in cases $M_3$, $M_5$, $N_1$ and $N_4$
is proved similarly to cases $M_2$, $M_4$, $N_2$ and $N_3$
correspondingly.
\end{rem}

Lemma \ref{mainlemma} is now proved.
\end{proof}
\vspace{0.2cm}

With this we now prove Theorem B. Observe that the trace functions
of isospectral compact Riemannian manifolds coincide
(Proposition \ref{pro1}) and the trace function $ \text{tr}(H_M) $
determines all parameters of the fundamental set of the corresponding
manifold $M$, consequently, it determines, up to isometry, the
manifold $M$ (Lemma \ref{mainlemma}).

Theorem B is now proved.

\section{The proof of Theorem C}

\begin{theoremC}
There is a unique family of pairs of isospectral non-homeomorphic flat $3$~-~manifolds
which consists of manifolds $M_4$ and $M_6$.
\end{theoremC}
\vspace{0.2cm}

\begin{proof}

In order to proof Theorem C it is sufficient to verify that
the trace functions $\text{tr}(H_{M_4})$ and $\text{tr}(H_{M_6})$
coincide under given values of the parameters of
$M_4$ and $M_6$, since the spectrum of the manifold $M$ is
exactly determined by the trace function $\text{tr}(H_M)$
(Proposition~\ref{pro2}).

Using Theorem A in case $M_4$ and $M_6$ we have

\begin{equation*}
\text{tr}(H_{M_4})
=
\frac{\text{vol}(M_4)}{(4\pi t)^{3/2}}
\sum\limits_{\boldsymbol{\ell} \in \Lambda}
e^{-|\boldsymbol{\ell}\,|^{2}/4t}
+
\frac{\ell_{1}}{4\sqrt{\pi t}}
\sum\limits_{m \in \mathbb{Z}}
e^{-(m + \frac{1}{4})^{2} \ell_{1}^{2} / 4t}
+
\frac{\ell_{1}}{8 \sqrt{\pi t}}
\sum\limits_{n \in \mathbb{Z}}
e^{-(n + \frac{1}{2})^{2} \ell_{1}^{2} / 4t},
\end{equation*}

\begin{multline*}
\text{tr}(H_{M_6})
=
\frac{\text{vol}(M_6)}{(4\pi
t)^{3/2}}\sum\limits_{\boldsymbol{b}\in\Lambda_1}e^{-|\boldsymbol{b}\,|^{2}/4t}
+
\frac{b_{1}}{8\sqrt{\pi t}}\sum\limits_{m\in
\mathbb{Z}}e^{-(m+\frac{1}{2})^{2}b_{1}^{2}/4t}
\\
+
\frac{b_{2}}{8\sqrt{\pi t}}\sum\limits_{n\in
\mathbb{Z}}e^{-(n+\frac{1}{2})^{2}b_{2}^{2}/4t}
+
\frac{b_{3}}{8\sqrt{\pi t}}\sum\limits_{k\in
\mathbb{Z}}e^{-(k+\frac{1}{2})^{2}b_{3}^{2}/4t}.
\end{multline*}

Here $\Lambda=
\{\boldsymbol{\ell}
=
m_1\boldsymbol{\ell}_{1} + m_2\boldsymbol{\ell}_{2}
+ m_3\boldsymbol{\ell}_{3},\ \ m_i\in \mathbb{Z}\,,\
\boldsymbol{\ell}_{i}\in \mathbb{R}^3,\, i=1,2,3\}$ is
a lattice for the torus $4$-fold covering $M_4$,
where
mutually orthogonal vectors
$\boldsymbol{\ell}_{1}/4$,
$\boldsymbol{\ell}_{2}$ and $\boldsymbol{\ell}_{3}$
form the fundamental set of the manifold $M_4$,
$\bell_1=(\ell_1,0,0)$.
Similarly,
$\Lambda_1=\{\boldsymbol{b}
=
m_1 \boldsymbol{b}_{1} + m_2 \boldsymbol{b}_{2}
+ m_3 \boldsymbol{b}_{3},\ \ m_i \in \mathbb{Z}\,,\
\boldsymbol{b}_{i}\in \mathbb{R}^3,\, i=1,2,3\}$ is a lattice
for the torus $2$-fold covering $M_2$,
which in turn $2$-fold covers the manifold $M_6$.
Vectors $\boldsymbol{b}_{1}/2$,
$\boldsymbol{b}_{2}/2$ and $\boldsymbol{b}_{3}$
form the fundamental set of $M_6$,
$\boldsymbol{b}_1=(b_1,0,0)$,
$\boldsymbol{b}_2=(0,b_2,0)$ and
$\boldsymbol{b}_3=(0,0,b_3)$.

Observe that if $b_2 = \ell_1$ then sums
\begin{equation*}
\frac{\ell_{1}}{8 \sqrt{\pi t}}
\sum\limits_{n \in \mathbb{Z}}
e^{-(n + \frac{1}{2})^{2} \ell_{1}^{2} / 4t}
\hspace{10pt}
\text{and}
\hspace{10pt}
\frac{b_{2}}{8\sqrt{\pi t}}\sum\limits_{n\in
\mathbb{Z}}e^{-(n+\frac{1}{2})^{2}b_{2}^{2}/4t}
\end{equation*}
coincide.

Also we notice that
\begin{multline*}
\sum\limits_{m \in \mathbb{Z}}
e^{-( m + \frac{1}{4} )^2 \ell_1^2 / 4 t}
=
\sum\limits_{m \in \mathbb{Z}}
e^{-( 2 m + \frac{1}{2} )^2 ( \ell_1 / 2 )^2 / 4 t}
=
\sum\limits_{m \in \mathbb{Z}}
e^{-( - 2 m + \frac{1}{2} )^2 ( \ell_1 / 2 )^2 / 4 t}
\\
=
\sum\limits_{m \in \mathbb{Z}}
e^{-( 2 m - \frac{1}{2} )^2 ( \ell_1 / 2 )^2 / 4 t}
=
\sum\limits_{m \in \mathbb{Z}}
e^{-( (2 m - 1) + \frac{1}{2} )^2 ( \ell_1 / 2 )^2 / 4 t}
=
\frac{1}{2}
\sum\limits_{k \in \mathbb{Z}}
e^{-( k + \frac{1}{2} )^2 (\ell_1 / 2)^2 / 4 t}.
\end{multline*}

Now it is easy to verify that if
$ b_1 = b_3 = \ell_1 / 2 $ then we have the following
equality.
\begin{equation*}
\frac{ \ell_{1} / 2 }{ 4 \sqrt{\pi t}}
\sum\limits_{k \in \mathbb{Z}}
e^{-( k + \frac{1}{2} )^{2} ( \ell_{1} / 2 )^{2} / 4t}
=
\frac{b_{1}}{8\sqrt{\pi t}}\sum\limits_{m\in
\mathbb{Z}}e^{-(m+\frac{1}{2})^{2}b_{1}^{2}/4t}
+
\frac{b_{3}}{8\sqrt{\pi t}}\sum\limits_{k\in
\mathbb{Z}}e^{-(k+\frac{1}{2})^{2}b_{3}^{2}/4t}.
\end{equation*}

We are left to check the coincidence of the following sums
\begin{equation*}
\frac{\text{vol}(M_4)}{(4\pi t)^{3/2}}
\sum\limits_{\boldsymbol{\ell} \in \Lambda}
e^{-|\boldsymbol{\ell}\,|^{2}/4t}
\hspace{10pt}
\text{and}
\hspace{10pt}
\frac{\text{vol}(M_6)}{(4\pi
t)^{3/2}}\sum\limits_{\boldsymbol{b}\in\Lambda_1}e^{-|\boldsymbol{b}\,|^{2}/4t}.
\end{equation*}

Since $\text{Vol}(M_4) = \text{Vol}(M_6)$, we have
$ b_1 = \ell $, where $\ell = |\boldsymbol{\ell_2}| = |\boldsymbol{\ell_3}|$
(Theorem \ref{orientable} in case of $M_4$). Here we have used that
$ \text{vol}(M_4)~=~\ell_1 \ell^2 /4 $, $ \text{vol}(M_6) = b_1 b_2 b_3 /4 $.
Consequently, $ \ell_1 = 2 \ell$. This gives us that lattices
$ \Lambda $ and $ \Lambda_1 $ are the same.

One-parameter family of pairs of isospectral non-homeomorphic compact flat
$3$-manifolds which consists of manifolds $M_4$ and $M_6$ is constructed.
Here the fundamental set of $M_4$ is formed by vectors
$ \boldsymbol{\ell}_1 / 4 $, $ \boldsymbol{\ell}_2 $ and $ \boldsymbol{\ell}_3 $,
where $ | \boldsymbol{\ell}_1 / 4 | = \ell / 2 $,
$ |\boldsymbol{\ell}_2| = |\boldsymbol{\ell}_3| = \ell $. The fundamental set of
$M_6$ is formed by vectors
$ \boldsymbol{b}_1 / 2 $, $ \boldsymbol{b}_2 / 2 $ and $ \boldsymbol{b}_3 $,
where $ | \boldsymbol{b}_1 / 2 | = \ell / 2 $,
$ |\boldsymbol{b}_2 / 2| = |\boldsymbol{b}_3| = \ell $.

Using the explicit formulas of $\text{tr}(H_M)$ of compact flat
$3$-manifolds~$M$ (Theorem A) it is easy to verify that
other examples of isospectral non-homeomorphic compact flat
$3$-manifolds do not exist.

Actually, it suffices to notice that the trace functions of non-orientable
manifolds contain terms with coefficient
$ \frac{1}{t} $ as distinct from the trace functions of orientable
manifolds which have the largest power with respect to
$t$ equals to $ \frac{1}{\sqrt{t}} $. Also observe that equality of
volumes of the manifolds and lattices of their covering tori
give us supplementary conditions.

Theorem C is now proved.

\end{proof}

\section{Appendix}

\subsection{The explicit formulas for the trace functions}

\setlength{\extrarowheight}{8pt}
\begin{center}
\begin{tabular}{|l|p{6in}|r|}
\multicolumn{2}{c}{Table 1}\\
\multicolumn{2}{c}{}\\
 \hline

$M_1$ & $F_1=\frac{\text{vol}(M_1)}{(4\pi t)^{3/2}}
\sum\limits_{\boldsymbol{\ell}\in\Lambda}e^{-|\boldsymbol{\ell}\,|^{2}/4t},$ \\
& $\Lambda$ is a lattice for the flat torus $M_1=\mathbb{R}^3/\Lambda$

\\ \hline

$M_2$ &  $F_2=\frac{\text{vol}(M_2)}{(4\pi
t)^{3/2}}\sum\limits_{\boldsymbol{\ell}\in\Lambda}e^{-|\boldsymbol{\ell}\,|^{2}/4t}
+
\frac{\ell_{1}}{4\sqrt{\pi t}}\sum\limits_{m\in
\mathbb{Z}}e^{-(m+\frac{1}{2})^{2}\ell_{1}^{2}/4t},$\\
& $\Lambda$
is a lattice for $M_1$, which is a two-fold
covering of the manifold $M_2$.
Here vectors $\boldsymbol{\ell}_{1}/2$,
$\boldsymbol{\ell}_{2}$ and $\boldsymbol{\ell}_{3}$
form the fundamental set of the manifold $M_2$,
$\bell_1=(\ell_1,0,0)$

\\  \hline

$M_3$ & $F_3=\frac{\text{vol}(M_3)}{(4\pi
t)^{3/2}}\sum\limits_{\boldsymbol{\ell}\in\Lambda}
e^{-|\boldsymbol{\ell}\,|^{2}/4t} +
\frac{\ell_{1}}{\sqrt{\pi t}}\sum\limits_{m\in
\mathbb{Z}}e^{-(m+\frac{1}{3})^{2}\ell_{1}^{2}/4t},$ \\
 & $\Lambda$
is a lattice for $M_1$, which is a three-fold
covering of the manifold $M_3$. Here vectors $\boldsymbol{\ell}_{1}/3$,
$\boldsymbol{\ell}_{2}$ and $\boldsymbol{\ell}_{3}$
form the fundamental set of the manifold $M_3$,
$\bell_1=(\ell_1,0,0)$

\\ \hline

$M_4$ & $F_4=\frac{\text{vol}(M_4)}{(4\pi t)^{3/2}}
\sum\limits_{\boldsymbol{\ell} \in \Lambda}
e^{-|\boldsymbol{\ell}\,|^{2}/4t}
+
\frac{\ell_{1}}{4\sqrt{\pi t}}
\sum\limits_{m \in \mathbb{Z}}
e^{-(m + \frac{1}{4})^{2} \ell_{1}^{2} / 4t}
+
\frac{\ell_{1}}{8 \sqrt{\pi t}}
\sum\limits_{n \in \mathbb{Z}}
e^{-(n + \frac{1}{2})^{2} \ell_{1}^{2} / 4t},$\\
& $\Lambda$
is a lattice for $M_1$, which
is a four-fold covering of the manifold $M_4$. Here vectors
$\boldsymbol{\ell}_{1}/4$, $\boldsymbol{\ell}_{2}$
and $\boldsymbol{\ell}_{3}$
form the fundamental set of the manifold $M_4$,
$\bell_1=(\ell_1,0,0)$

\\ \hline

$M_5$ & $F_5 = \frac{\text{vol}(M_5)}{(4\pi t)^{3/2}}
\sum\limits_{\boldsymbol{\ell} \in \Lambda}
e^{-|\boldsymbol{\ell}\,|^{2}/4t}
+
\frac{\ell_{1}}{2\sqrt{\pi t}}
\sum\limits_{m \in \mathbb{Z}}
e^{-(m + \frac{1}{6})^{2} \ell_{1}^{2} / 4t} $\\
&
$
+
\frac{\ell_{1}}{2 \sqrt{\pi t}}
\sum\limits_{n \in \mathbb{Z}}
e^{-( n + \frac{1}{3})^{2} \ell_{1}^{2} / 4t}
+
\frac{\ell_{1}}{4 \sqrt{\pi t}}
\sum\limits_{k \in \mathbb{Z}}
e^{-( k + \frac{1}{2})^{2} \ell_{1}^{2} / 4t},$ \\
& $\Lambda$ is a lattice for $M_1$, which is a six-fold covering of the
manifold $M_5$. Here vectors $\boldsymbol{\ell}_{1}/6$,
$\boldsymbol{\ell}_{2}$ and
$\boldsymbol{\ell}_{3}$
form the fundamental set of the manifold $M_5$,
$\bell_1=(\ell_1,0,0)$

\\ \hline

\end{tabular}
\end{center}
\setlength{\extrarowheight}{8pt}
\begin{center}
\begin{tabular}{|l|p{6in}|r|}
\multicolumn{2}{c}{Table 1: continue}\\
\multicolumn{2}{c}{} \\
\hline

$M_6$ & $ F_6=\frac{\text{vol}(M_6)}{(4\pi
t)^{3/2}}\sum\limits_{\boldsymbol{\ell}\in\Lambda}e^{-|\boldsymbol{\ell}\,|^{2}/4t}
+
\frac{\ell_{1}}{8\sqrt{\pi t}}\sum\limits_{m\in
\mathbb{Z}}e^{-(m+\frac{1}{2})^{2}\ell_{1}^{2}/4t}$
\\
& $+
\frac{\ell_{2}}{8\sqrt{\pi t}}\sum\limits_{n\in
\mathbb{Z}}e^{-(n+\frac{1}{2})^{2}\ell_{2}^{2}/4t}+
\frac{\ell_{3}}{8\sqrt{\pi t}}\sum\limits_{k\in
\mathbb{Z}}e^{-(k+\frac{1}{2})^{2}\ell_{3}^{2}/4t},$ \\
& $\Lambda$ is a lattice for $M_1$, which is a two-fold covering of the
manifold $M_2$, which in turn is a two-fold covering of the
manifold $M_6$. Here vectors $\boldsymbol{\ell}_{1}/2$,
$\boldsymbol{\ell}_{2}/2$ and
$\boldsymbol{\ell}_{3}$
form the fundamental set of the manifold $M_6$,
$\bell_1=(\ell_1,0,0)$,
$\bell_2=(0,\ell_2,0)$ and
$\bell_3=(0,0,\ell_3)$

\\ \hline

$N_1$ & $P_1=\frac{\text{vol}(N_1)}{(4\pi
t)^{3/2}}\sum\limits_{\boldsymbol{\ell}\in\Lambda}e^{-|\boldsymbol{\ell}\,|^{2}/4t}
+
\frac{\ell_{1}\ell_{2}\sin\varphi}{8\pi t}\sum\limits_{m,\, n\in
\mathbb{Z}}e^{-[(\ell_1(m+\frac{1}{2})+\ell_2n\cos\varphi)^{2}+(n\ell_{2}\sin\varphi)^{2}]/4t},$\\

& $\Lambda$ is a lattice for the flat torus
$M_1=\mathbb{R}^3/\Lambda$, which is a two-fold covering of the
manifold $N_1$. Here vectors $\boldsymbol{\ell}_{1}/2$,
$\boldsymbol{\ell}_{2}$ and
$\boldsymbol{\ell}_{3}$
form the fundamental set of the manifold $N_1$,
$\boldsymbol{\ell}_{1}=(\ell_1,0,0)$,
$\boldsymbol{\ell}_{2}=(\ell_2\cos\varphi,\ell_2\sin\varphi,0)$
and
$\boldsymbol{\ell}_{3}=(0,0,\ell_3)$, where
$\varphi$ is an angle between $\boldsymbol{\ell}_{1}$ and
$\boldsymbol{\ell}_{2}$, $0 < \varphi < \pi$

\\ \hline

$N_2$ & $P_2=\frac{\text{vol}(N_2)}{(4\pi
t)^{3/2}}\sum\limits_{\boldsymbol{\ell}\in\Lambda}e^{-|\boldsymbol{\ell}\,|^{2}/4t}
+
\frac{ \ell_{1} \ell_{2} \sin\varphi }{ 16 \pi t}
\sum\limits_{ m,\,n \in \mathbb{Z}}
e^{-[( \ell_1 ( m + \frac{1}{2} ) + n \ell_2 \cos\varphi)^2
+
( n \ell_2 \sin\varphi)^2 / 4t ]}$
\\
&$+
\frac{ \ell_{1} \ell_{2} \sin\varphi }{ 16 \pi t}
\sum\limits_{ m,\,n \in \mathbb{Z}}
e^{-[ \ell_1  m  + ( n + \frac{1}{2} ) \ell_2 \cos\varphi)^2
+
( ( n + \frac{1}{2} ) \ell_2 \sin\varphi)^2 / 4t ]}$
\\
& $\Lambda$ is a lattice
for $M_1$, which is a two-fold
covering of the manifold $N_2$. Here
vectors $\boldsymbol{\ell}_{1}/2$, $\boldsymbol{\ell}_{2}$ and
$\boldsymbol{\ell}_{3}$
form the fundamental set of the manifold $N_2$,
$\boldsymbol{\ell}_{1}=(\ell_1,0,0)$,

\noindent
$\boldsymbol{\ell}_{2}=(\ell_2\cos\varphi,\ell_2\sin\varphi,0)$
and $\boldsymbol{\ell}_{3}=(\frac{1}{2}(\ell_1+\ell_2\cos\varphi),
\frac{1}{2}\ell_2\sin\varphi,h)$, where $\varphi$ is
an angle between $\boldsymbol{\ell}_{1}$ and
$\boldsymbol{\ell}_{2}$,
$0 < \varphi < \pi$, $h$ is a projection of $\boldsymbol{\ell_3}$ on
$\boldsymbol{e_3}$

\\ \hline

$N_3$ & $P_3=\frac{\text{vol}(N_3)}{(4\pi
t)^{3/2}}\sum\limits_{\boldsymbol{\ell}\in\Lambda}e^{-|\boldsymbol{\ell}\,|^{2}/4t}
+
\frac{\ell_{1}}{8\sqrt{\pi t}}\sum\limits_{m\in
\mathbb{Z}}e^{-[(m+\frac{1}{2})^{2}\ell_{1}^{2}]/4t}$
\\
& $+ \frac{\ell_1\ell_{2}}{16\pi t}\sum\limits_{m,n\in
\mathbb{Z}}e^{-[\ell_1^{2}m^2+\ell_{2}^{2}(n+\frac{1}{2})^2]/4t}+
\frac{\ell_1\ell_{3}}{16\pi t}\sum\limits_{m,k\in
\mathbb{Z}}e^{-[\ell_1^2(m+\frac{1}{2})^{2}+\ell_{3}^{2}k^2]/4t},$\\
& $\Lambda$ is a lattice for
$M_1$, which is a two-fold covering of the
manifold $M_2$, which in turn is a two-fold covering of the
manifold $N_3$. Here vectors $\boldsymbol{\ell}_{1}/2$,
$\boldsymbol{\ell}_{2}/2$ and
$\boldsymbol{\ell}_{3}$
form the fundamental set of the manifold $N_3$,
$\bell_1=(\ell_1,0,0)$,
$\bell_2=(0,\ell_2,0)$ and
$\bell_3=(0,0,\ell_3)$

\\ \hline

$N_4$ & $P_4=\frac{\text{vol}(N_4)}{(4\pi
t)^{3/2}}\sum\limits_{\boldsymbol{\ell}\in\Lambda}e^{-|\boldsymbol{\ell}\,|^{2}/4t}
+
\frac{\ell_{1}}{8\sqrt{\pi t}}\sum\limits_{m\in
\mathbb{Z}}e^{-[(m+\frac{1}{2})^{2}\ell_{1}^{2}]/4t}$
\\
& $+
\frac{\ell_1\ell_{2}}{16\pi t}\sum\limits_{m,n\in
\mathbb{Z}}e^{-[\ell_1^{2}m^2+\ell_{2}^{2}(n+\frac{1}{2})^2]/4t}+
\frac{\ell_1\ell_{3}}{16\pi t}\sum\limits_{m,k\in
\mathbb{Z}}e^{-[\ell_1^2(m+\frac{1}{2})^{2}+\ell_{3}^{2}(k+\frac{1}{2})^2]/4t},$\\
& $\Lambda$ is a lattice for
$M_1$, which is a two-fold covering of the
manifold $M_2$, which in turn is a two-fold covering of the
manifold $N_4$. Here vectors $\boldsymbol{\ell}_{1}/2$,
$\boldsymbol{\ell}_{2}/2$ and
$\boldsymbol{\ell}_{3}$
form the fundamental set of the manifold $N_4$,
$\bell_1=(\ell_1,0,0)$,
$\bell_2=(0,\ell_2,0)$ and
$\bell_3=(0,0,\ell_3)$

\\ \hline

\end{tabular}
\end{center}

\newpage

\subsection{The fundamental sets of flat 3-manifolds}

Here we give figures of the fundamental sets of flat
$3$-manifolds. In order to obtain a ma\-ni\-fold~$M$ one must glue
faces of the corresponding fundamental set~$M^{\ast}$
which are labelled with the same letter.

\begin{figure}
[h]
\centering
\includegraphics
[width=7cm, height=3cm]
{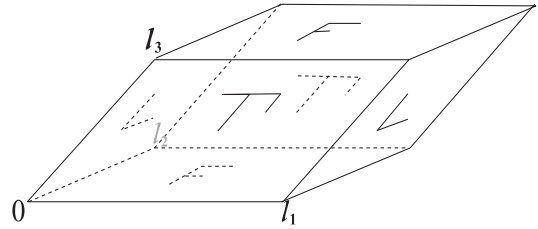}
\caption{The fundamental set of $M_1$}
\end{figure}

\begin{figure}[h]
\includegraphics[width=15cm, height=3cm]{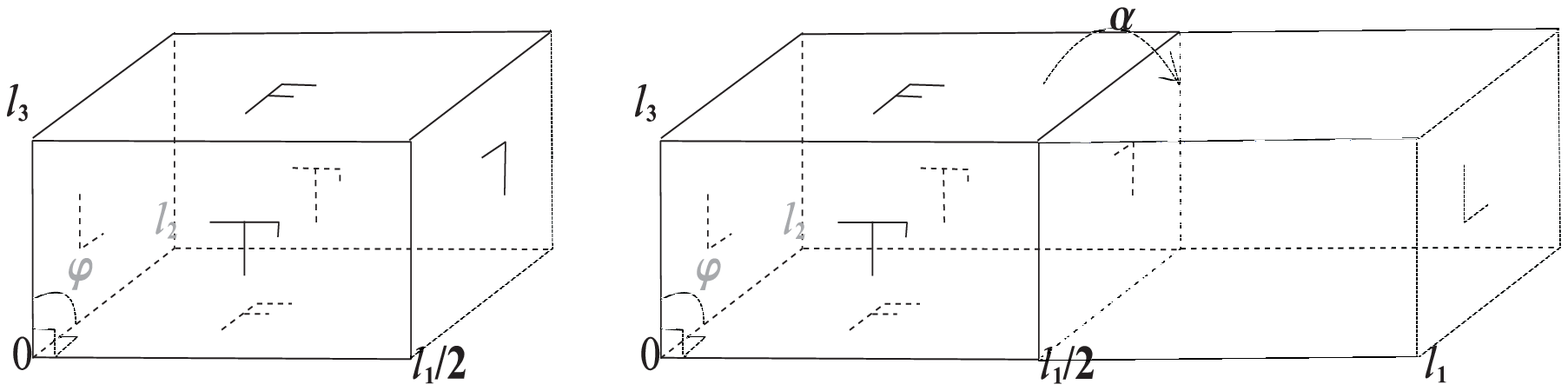}
\caption{The fundamental set of $M_2$ and covering torus $M_1$,
$M_1 \stackrel{2}{\rightarrow}M_2$}
\end{figure}

\begin{figure}[h]
\includegraphics[width=15cm, height=5cm]{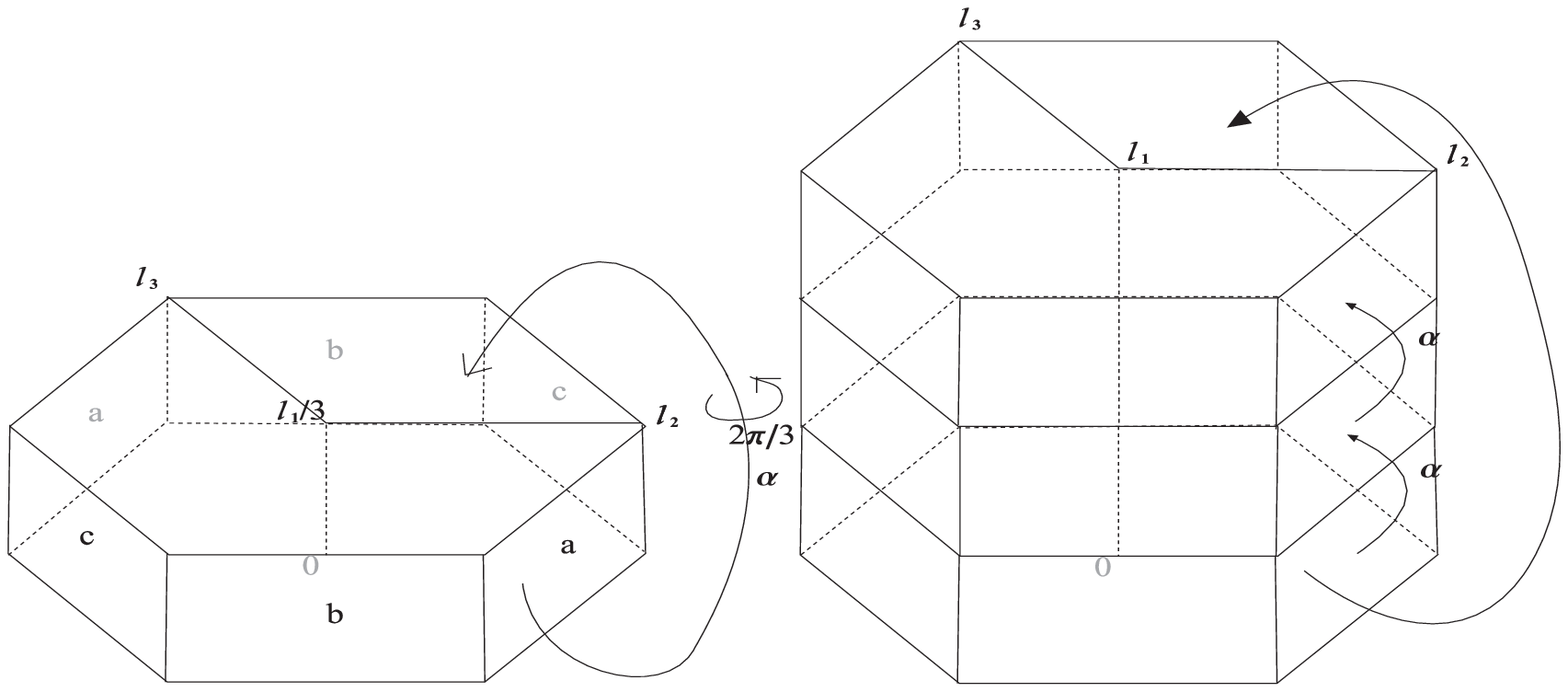}
\caption{The fundamental set of $M_3$ and covering torus $M_1$,
$M_1 \stackrel{3}{\rightarrow}M_3$}
\end{figure}

\begin{figure}
\includegraphics[width=15cm, height=3cm]{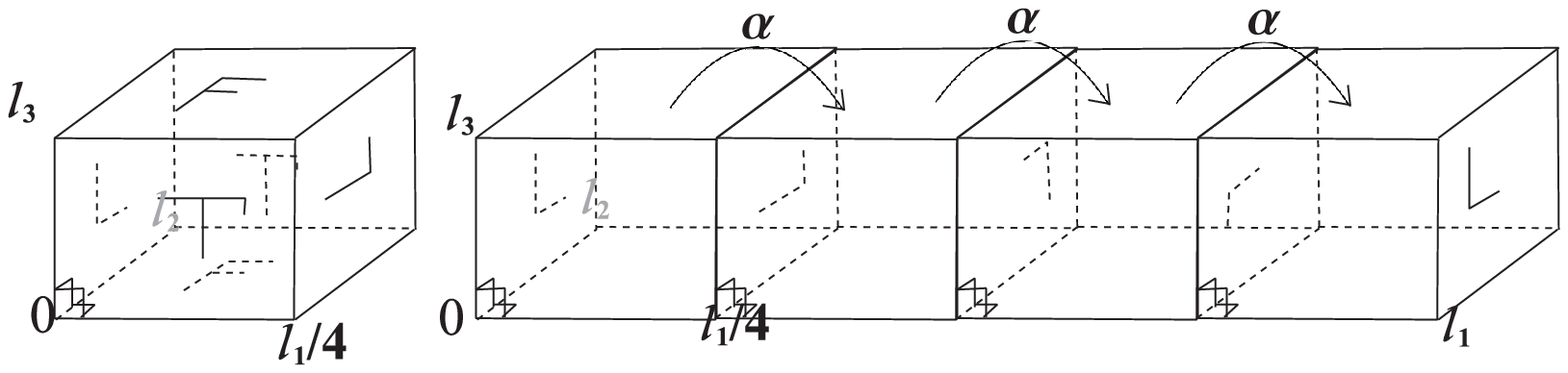}
\caption{The fundamental set of $M_4$ and covering torus $M_1$,
$M_1 \stackrel{4}{\rightarrow}M_4$}
\end{figure}

\begin{figure}
\includegraphics[width=15cm, height=5cm]{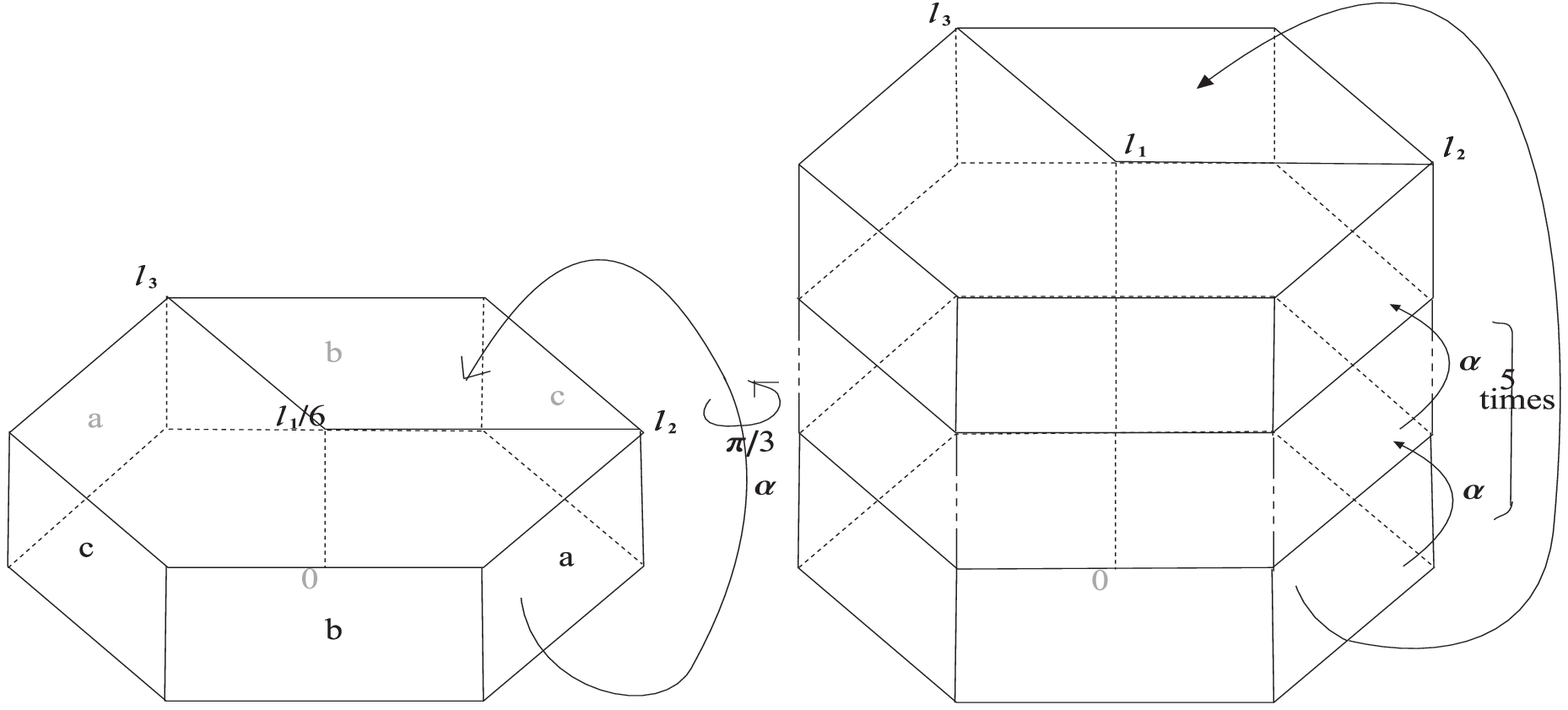}
\caption{The fundamental set of $M_5$ and covering torus $M_1$,
$M_1 \stackrel{6}{\rightarrow}M_5$}
\end{figure}

\begin{figure}
\includegraphics[width=15cm, height=5cm]{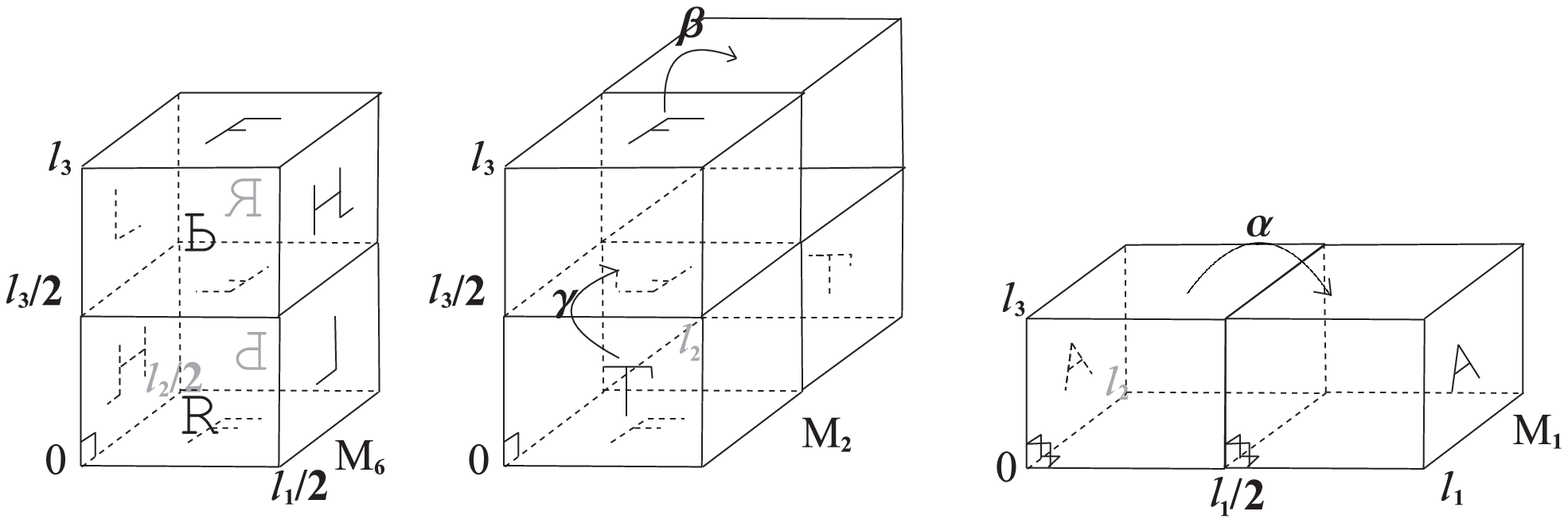}
\caption{The fundamental set of $M_6$ and covering manifolds
$M_1 \stackrel{2}{\rightarrow}M_2\stackrel{2}{\rightarrow}M_6$}
\end{figure}

\begin{figure}
\includegraphics[width=15cm, height=3cm]{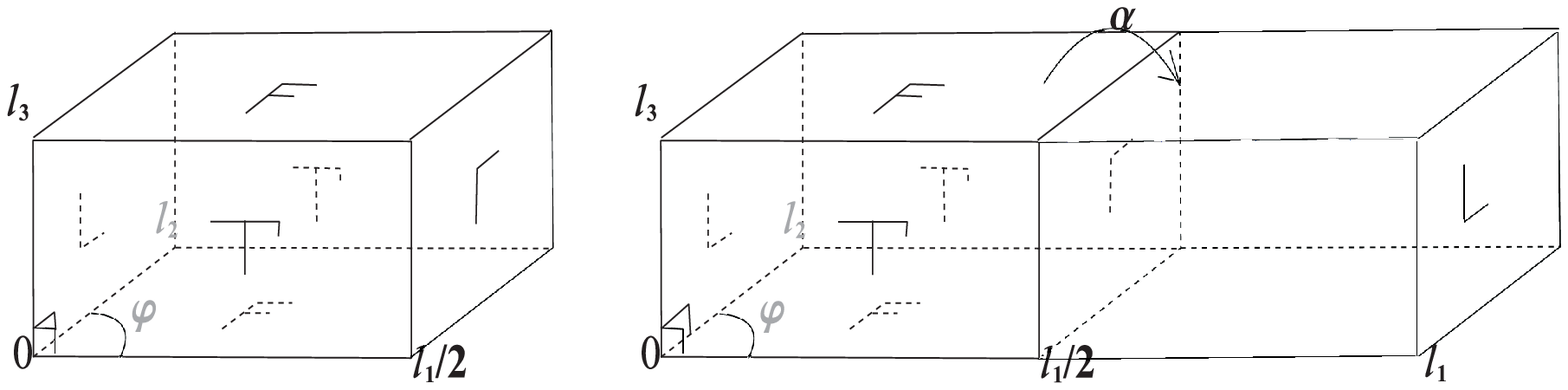}
\caption{The fundamental set of $N_1$ and covering torus $M_1$,
$M_1 \stackrel{2}{\rightarrow}N_1$}
\end{figure}

\begin{figure}
\includegraphics[width=15cm, height=3cm]{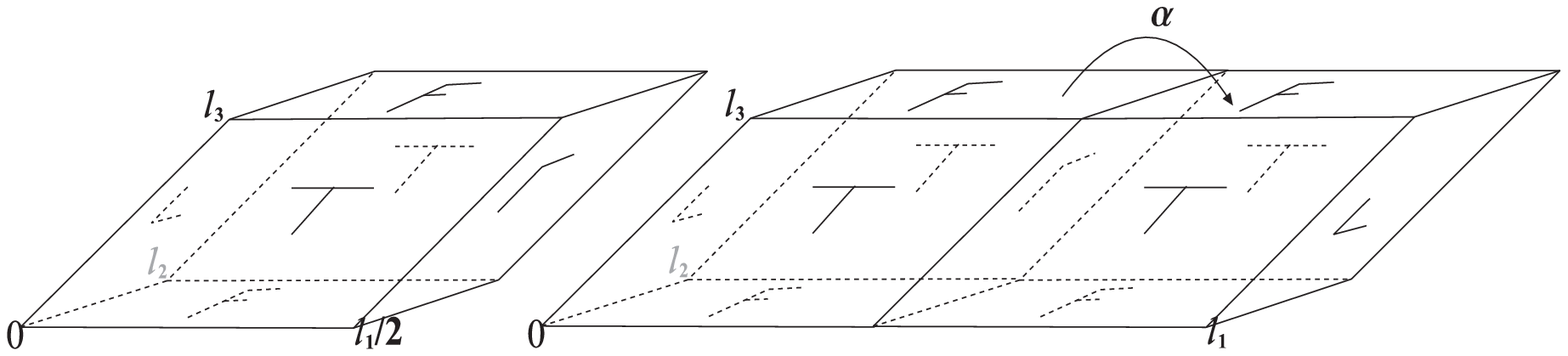}
\caption{The fundamental set of $N_2$ and covering torus $M_1$,
$M_1 \stackrel{2}{\rightarrow}N_2$}
\end{figure}

\begin{figure}
\includegraphics[width=15cm, height=5cm]{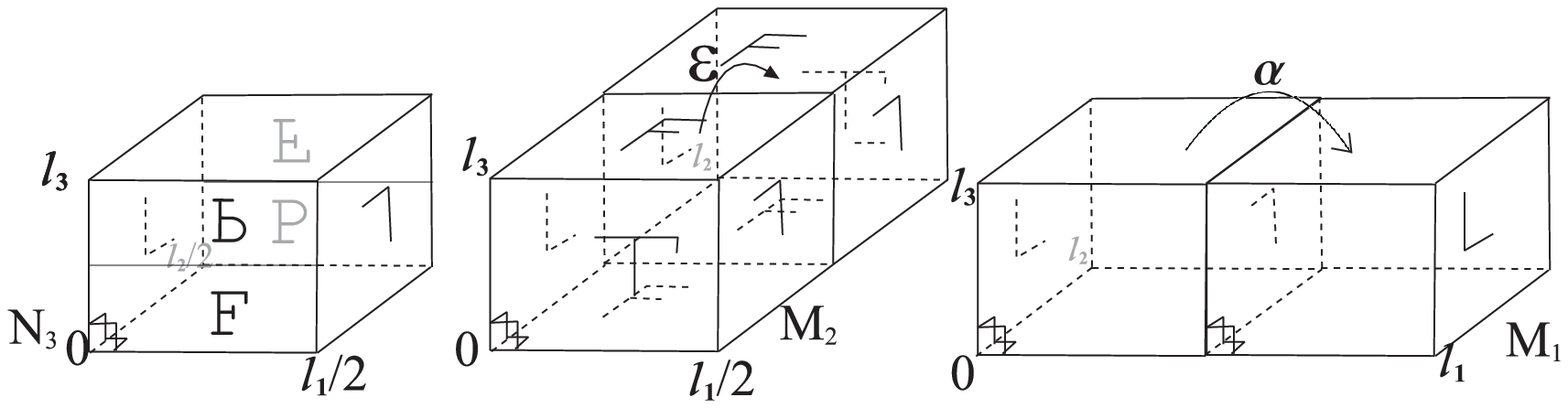}
\caption{The fundamental set of $N_3$ and covering manifolds
$M_1 \stackrel{2}{\rightarrow}M_2\stackrel{2}{\rightarrow}N_3$}

\end{figure}

\begin{figure}
\includegraphics[width=15cm, height=5cm]{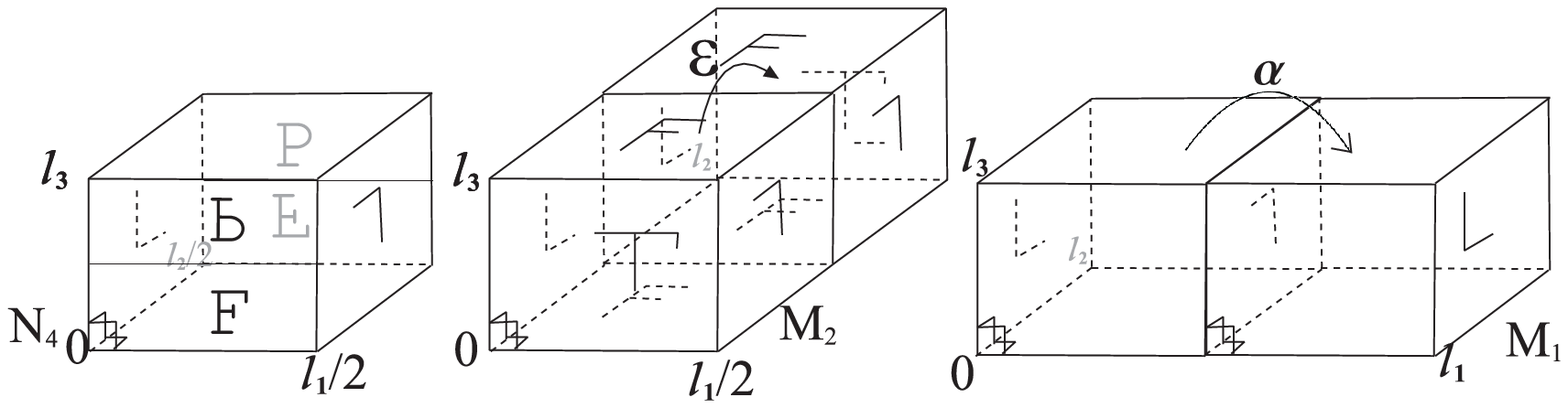}
\caption{The fundamental set of $N_4$ and covering manifolds
$M_1 \stackrel{2}{\rightarrow}M_2\stackrel{2}{\rightarrow}N_4$}

\end{figure}

\newpage

\addcontentsline{toc}{section}{\hspace{15pt} References}

\end{document}